\documentclass[final]{siamltex}
\usepackage{epsfig}
\usepackage{amssymb}

\newtheorem{example}[theorem]{Example}
\newtheorem{problem}[theorem]{Problem}
\def\IC{{\mathbb C}}
\def\IR{{\mathbb R}}
\def\IN{{\mathbb N}}
\def\cL{{\mathcal L}}
\def\cU{{\mathcal U}}
\def\cV{{\mathcal V}}

\def\cH{{\mathcal H}}

\def\cH{{\mathcal H}}
\def\cP{{\mathcal P}}

\def\cS{{\mathcal S}}

\def\cT{{\mathcal T}}
\def\cW{{\mathcal W}}

\def\({\left(}
\def\){\right)}
\def\diag{{\rm diag}\,}
\def\conv{{\rm conv}\,}

\def\Re{{\rm Re}\,}
\def\Im{{\rm Im}\,}

\def\int{{\bf Int}\,}

\title {Higher rank numerical ranges of normal matrices}
\author{Hwa-Long Gau\thanks{
Department of Mathematics, National Central University,
Chung-Li 320, Taiwan
(Email: hlgau@math.ncu.edu.tw).
Research of Gau was supported by National Science Council of the Republic
of China.}
\and
Chi-Kwong Li\thanks{
Department of Mathematics, College of William \& Mary,
Williamsburg, VA 23185
(Email: ckli@math.wm.edu).
Research of Li was supported by USA NSF.
Li was also supported by the William and Mary Plumeri Award.
Li is an honorary professor of the University of Hong Kong
and an honorary professor of the Taiyuan University of Technology.}
\and
Yiu-Tung Poon\thanks{
Department of Mathematics, Iowa State University,
Ames, IA 50051
(Email: ytpoon@iastate.edu).
Research of Poon was supported by USA NSF.}
\and
Nung-Sing Sze\thanks{
Department of Applied Mathematics,
The Hong Kong Polytechnic University,
Hung Hom, Hong Kong
(Email: raymond.sze@inet.polyu.edu.hk).
Research of Sze was supported by the Hong Kong Polytechnic University startup grant.}
}

\begin{document}
\maketitle

%\date{}

\begin{abstract}
The higher rank numerical range is closely connected to  the construction of
quantum error correction code for a noisy quantum channel. It is
known that if a normal matrix $A \in M_n$ has eigenvalues $a_1,
\dots, a_n$, then its higher rank numerical range $\Lambda_k(A)$ is
the intersection of convex polygons with
vertices $a_{j_1}, \dots, a_{j_{n-k+1}}$, where $1 \le j_1 <
\cdots < j_{n-k+1} \le n$. In this paper, it is shown that the
higher rank numerical range of a normal matrix with $m$ distinct
eigenvalues can be written as the intersection of no more than
$\max\{m,4\}$ closed half planes.
In addition, given a convex polygon $\cP$ a
construction is given for a normal matrix $A \in M_n$ with
minimum $n$ such that $\Lambda_k(A) = \cP$. In particular, if
$\cP$ has $p$ vertices, with $p \ge 3$, there is a normal matrix
$A \in M_n$ with $n \le \max\left\{p+k-1,  2k+2 \right\}$ such
that $\Lambda_k(A) = \cP$.
\end{abstract}

\begin{AMS}
{15A60, 15A90, 47N50, 81P68}
\end{AMS}


\section{Introduction}
\setcounter{equation}{0}

Let $M_n$ be the algebra of $n\times n$ complex matrices regarded as
linear operators acting on the $n$-dimensional Hilbert space $\IC^n$.
The {\em classical numerical range} of $A\in M_n$ is defined and denoted by
$$W(A) = \{ x^*Ax\in \IC:
x\in \IC^n \hbox{ with } x^*x = 1\},$$ which is a  useful concept
in studying matrices and operators; see \cite{HJ}.

In the context of quantum information theory,
if the quantum states are represented as matrices in $M_n$,
then a {\it quantum channel} is
a trace preserving completely positive map
$L: M_n \rightarrow M_n$ with the following
operator sum representation
\begin{equation} \label{opersum}
L(A) = \sum_{j=1}^r E_j^* A E_j,
\end{equation}
where $E_1, \dots, E_r \in M_n$ satisfy
$\sum_{j=1}^r E_j E_j^* = I_n$.
The matrices $E_1, \dots, E_r$ are known as the
{\em error operators} of the quantum channel $L$.
A subspace $V$ of $\IC^n$
is a {\em quantum error correction code} for the channel $L$
if and only if the orthogonal projection $P \in M_n$
with range space $V$ satisfies $PE_i^*E_jP = \gamma_{ij} P$
for all $i,j \in \{1, \dots, r\}$; for example, see
\cite{KL,KLV,KLPL}.
In this connection, for $1 \le k < n$ researchers
define the {\em rank-$k$ numerical range} of $A\in M_n$ by
$$\Lambda_k(A) = \{ \lambda \in \IC: PAP = \lambda P
\hbox{ for some rank-$k$ orthogonal  projection } P\},$$
and the {\em joint rank-$k$ numerical range} of $A_1, \dots, A_m
\in M_n$ by $\Lambda_k(A_1, \dots, A_m)$ to be the collection of
complex vectors $(a_1, \dots, a_m) \in \IC^{1\times m}$ such that
$PA_jP = a_j P$ for a rank-$k$ orthogonal projection $P \in M_n$.
Evidently, there is a quantum error correction code $V$ of dimension
$k$ for the quantum channel $L$ described in (\ref{opersum})
 if and only if
$\Lambda_k(A_1, \dots, A_m)$ is non-empty for $(A_1, \dots, A_m)
= (E_1^*E_1, E_1^*E_2,  \dots, E_r^*E_r)$. Also, it is easy to see
that if $(a_1, \dots, a_m) \in \Lambda_k(A_1, \dots, A_m)$ then
$a_j \in \Lambda_k(A_j)$ for $j = 1,\dots, m$.
When $k = 1$,  $\Lambda_k(A)$
reduces to the classical numerical range $W(A)$.

Recently,
interesting results have been obtained for the rank-$k$ numerical
range and the joint rank-$k$ numerical range; see
\cite{Cet,Cet0,Cet1,Cet2,GLW,LP,LPS,LPS2,LS,W1}. In particular,
an explicit description of the rank-$k$ numerical range of $A \in
M_n$ is given in \cite{LS}, namely,
\begin{eqnarray}\label{eq1.1}
\Lambda_k(A) = \bigcap_{\xi\in [0,2\pi)} \{ \mu \in \IC:
e^{-i\xi}\mu+ e^{i\xi}\overline{\mu} \le
\lambda_k(e^{-i\xi}A + e^{i\xi}A^*) \},
\end{eqnarray}
where $\lambda_k(X)$ is the $k$th largest eigenvalue of
a Hermitian matrix $X$.

In the study of quantum error correction, there are channels such
as the randomized unitary channels and Pauli channels whose error
operators are commuting normal matrices. Thus, it is of interest
to study the rank-$k$ numerical ranges of normal matrices.
Although   the error operators of a generic quantum channel
may   not commute, a good understanding
of the special case would lead to deeper insights and
more proof techniques for the general case.

Given $S\subseteq \IC$, let $\conv S$ denote the smallest convex
subset of $\IC$ containing $S$. For a normal matrix $A \in M_n$
with eigenvalues $a_1, \dots, a_n$, it was conjectured in
\cite{Cet1,Cet2} that
\begin{equation}\label{eq1.2}
\Lambda_k(A)
= \bigcap_{1 \le j_1 < \cdots < j_{n-k+1} \le n}
\conv\{a_{j_1},\dots,a_{j_{n-k+1}}\},
\end{equation}
which is a convex polygon including its interior (if it is non-empty).
This conjecture was confirmed in \cite{LS} using the description
of $\Lambda_k(A)$ in (\ref{eq1.1}).
{\it In our discussion,
a polygon would always mean a convex polygon
with its interior.}

In this paper, we improve the description (\ref{eq1.2}) of the
rank-$k$ numerical range of a normal matrix. In particular,
in Section 2 we show that for a normal matrix $A$ with $m$
distinct eigenvalues, $\Lambda_k(A)$ can be written as the
intersection of no more than $\max\{m,4\}$ closed half planes in $\IC$.
Moreover, if $\Lambda_k(A) \ne \emptyset$, then it is a polygon
with no more than $m$ vertices.
We then consider the ``inverse'' problem, namely, for a given
polygon $\cP$, construct a normal matrix $A \in M_n$ with
$\Lambda_k(A) = \cP$. In other words, we study the necessary
condition for the existence of quantum channels whose error
operators have prescribed rank-$k$ numerical ranges. It is easy
to check that $\Lambda_k(\tilde A) = \cP$ if $\tilde A = A
\otimes I_k$ with $W(A) = \cP$.
Our goal is to find a normal matrix $\hat A$ with smallest size
so that $\Lambda_k(\hat A) = \cP$. To achieve this, we give a
necessary and sufficient condition for the existence of a normal
matrix $A\in M_n$ so that $\Lambda_k(A) = \cP$ in terms of
$k$-regular sets in $\IC$ (see Definition \ref{3.1}).
Furthermore, we show that the problem of finding a desired normal
matrix $A$ is equivalent to a combinatorial problem of extending a given
$p$ element set of unimodular complex numbers to a $k$-regular
set. We then give the solution of the problem in Section 4. As a
consequence of our results, if $\cP$ is a  polygon with $p$
vertices, then there is a normal matrix $A \in M_n$ with $$n
\le \max\left\{p+k-1,  2k+2 \right\}$$ such that $\Lambda_k(A) =
\cP$. Moreover, this upper bound is best possible in the sense
that there exists $\cP$ so that there is no matrix of smaller
dimension with rank-$k$ numerical range equal to $\cP$.

\section{Construction of higher rank numerical ranges}

By (\ref{eq1.1}),  $\Lambda_k(A)$ can be obtained as the intersection of
infinitely many closed half planes for a given $A \in M_n$.
Suppose $A$ is normal. By (\ref{eq1.2}), one can write $\Lambda_k(A)$
as the intersection of ${n\choose k-1}$ convex polygons so that
$\Lambda_k(A)$ is a polygon.  In particular,
it is well known that $\Lambda_1(A) = \conv\{a_1, \dots, a_m\}$,
where $a_1,\dots,a_m$ are the distinct eigenvalues of $A$.

There are nice interplay between the
algebraic properties of $A \in M_n$
and the geometric properties of $\Lambda_1(A) = W(A)$.
For instance,  $\Lambda_1(A)$
is always non-empty; $\Lambda_1(A)$ is a singleton if and only if
$A$ is a scalar matrix; $\Lambda_1(A)$ is a non-degenerate line
segment if
and only if $A$ is a non-scalar normal
matrix and its eigenvalues lie on a straight line.
Unfortunately, these results have no analogs for $\Lambda_k(A)$ if $k>1$.
First, the set $\Lambda_k(A)$ may be empty, see \cite{LPS2};
there are non-scalar matrices $A$ such that $\Lambda_k(A)$ is a singleton;
and there are non-normal matrices $A$ such that
$\Lambda_k(A)$ is a line segment.
Even for a normal matrix $A$,
it is not easy to determine whether $\Lambda_k(A)$ is empty,
a point or a line segment without actually constructing the
set $\Lambda_k(A)$. Moreover, there is no easy way to express
the vertices of the polygon $\Lambda_k(A)$ (if it is non-empty)
in terms of the eigenvalues of the normal matrix $A$ as in the case
of $\Lambda_1(A)$. Of course, one can use
(\ref{eq1.2}) to construct $\Lambda_k(A)$ for the normal
matrix $A$, but the number of polygons needed
in the construction will grow exponentially for large $n$ and $k$.
In the following, we will study efficient ways to generate
$\Lambda_k(A)$ for a normal matrix $A \in M_n$.
While it is difficult to use the eigenvalues of $A$ to
determine the set $\Lambda_k(A)$, it turns out that
we can use half planes determined by
the eigenvalues to generate $\Lambda_k(A)$ efficiently.
In the following, we will focus on the following problem.

\smallskip
\begin{problem} \rm
Determine the minimum number of half planes needed to construct
$\Lambda_k(A)$ using the eigenvalues of the normal matrix $A \in M_n$.
\end{problem}

\smallskip
As by-products, we will show that for a normal matrix $A$ with $m$ distinct
eigenvalues, $\Lambda_k(A)$ is either empty or is a polygon with
at most $m$ vertices.   In fact, by
examining the location of the eigenvalues of $A$ on the complex
plane, one may further reduce the number of half planes needed to
construct $\Lambda_k(A)$.

%\medskip
Suppose the eigenvalues of $A\in M_n$ are collinear. Then by a
translation, followed by a rotation, we may assume that $A$ is
Hermitian with eigenvalues  $a_1\ge \dots \ge a_n$. Then we have
$\Lambda_k(A)=\left[a_{n-k+1},a_{k}\right]$. So we focus on those normal matrices  whose eigenvalues are not collinear.

Let us motivate our result with the following examples,
which can be verified by using (\ref{eq1.2}).

\smallskip
\begin{example} \label{ex1} \rm
Let $A = \diag(1, w, w^2, \dots, w^{n-1})$ with
$w=e^{2\pi i/n}$. Then for $k \le n/2$,
we have $\Lambda_k(A) = \cap_{j=0}^{n-1} \cH_j$,
where
$$\cH_j = \left\{z\in \IC: \Re \left( e^{-\frac{(2j+k)\pi i}{n}} z \right) \le \cos \frac{k \pi}{n} \right\},$$
and only a small part of
$\conv\{w^{j-1}, w^{j-1+k}\}$ lies in $\Lambda_k(A)$.
\end{example}

\begin{center}
\epsfig{file=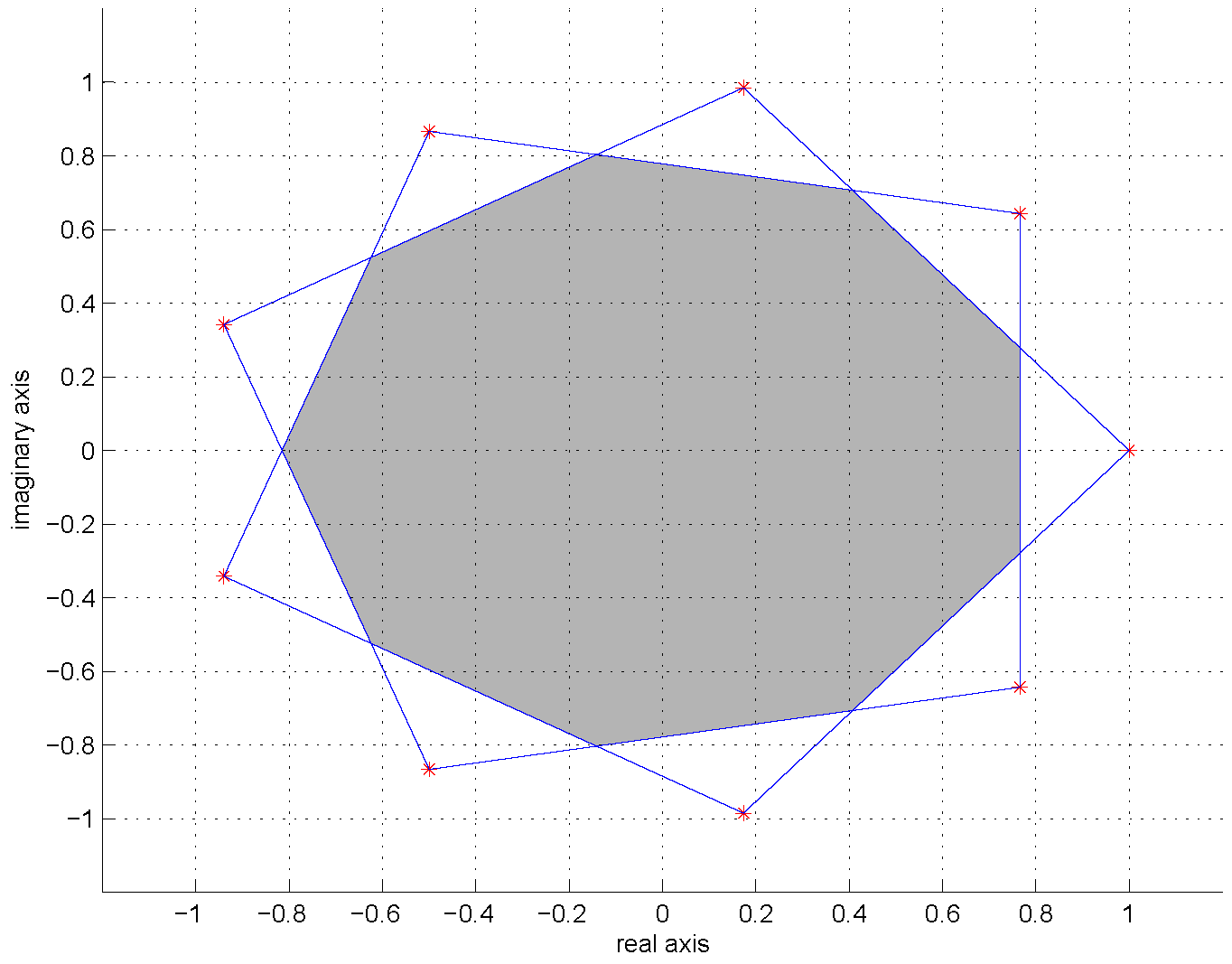, width =2in, height=2in} \qquad
\epsfig{file=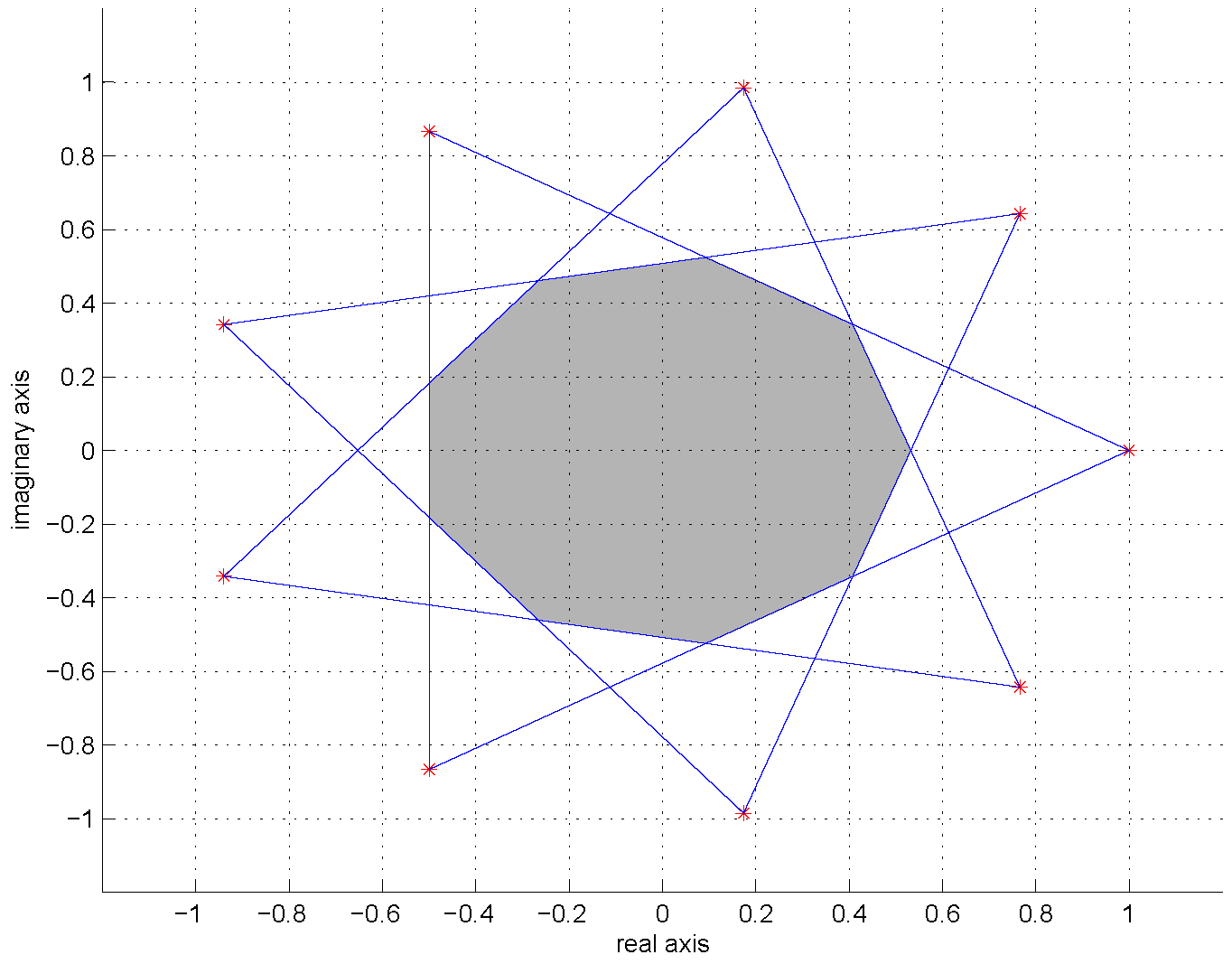, width =2in, height=2in}

$\Lambda_2(A)$ with $n = 9$ in Example \ref{ex1}\quad
$\Lambda_3(A)$ with $n = 9$ in Example \ref{ex1}
\end{center}

\smallskip
More generally, we have the following.

\smallskip
\begin{example}\label{eg1} \rm
Let $a_1,\dots,a_n$ be the eigenvalues of $A \in M_n$, with $n\ge 3$.
Suppose $\conv \{a_1,\dots,a_n\} = \cP$ is an $n$-sided convex polygon
containing the origin in the interior.
We may assume that $a_1,\dots,a_n$ are arranged in
the counterclockwise direction on the boundary of $\cP$.
For $j\in \{1, \dots, n\}$,
let $L_j$ be the line passing through $a_j$ and $a_{j+k}$,
where $a_{j+k} = a_{j+k-n}$ if $j+k> n$,
and $\cH_j$ be the closed half plane determined by $L_j$
which does not contain $a_\ell$ for $j< \ell <j+k$.
Then
$$\Lambda_k(A)=\bigcap_{j=1}^n\cH_j.$$
\end{example}

Note that each $\cH_j$ in Example \ref{eg1} contains exactly $n-k+1$
eigenvalues of $A$.

The situation is more complicated if $\Lambda_1(A)$ is not an $n$-sided
convex polygon for the normal matrix $A \in M_n$.

\smallskip
\begin{example} \label{ex2} \rm
Suppose $B = \diag(1,i,-1,-i, 2,2i,-2,-2i,3,3i,-3,-3i)$.
One can see from the figures that
the eigenvalues $1,i,-1,-i$ are interior points of $\Lambda_2(B)$
while these eigenvalues are the vertices of $\Lambda_3(B)$.
\end{example}

\begin{center}
\epsfig{file=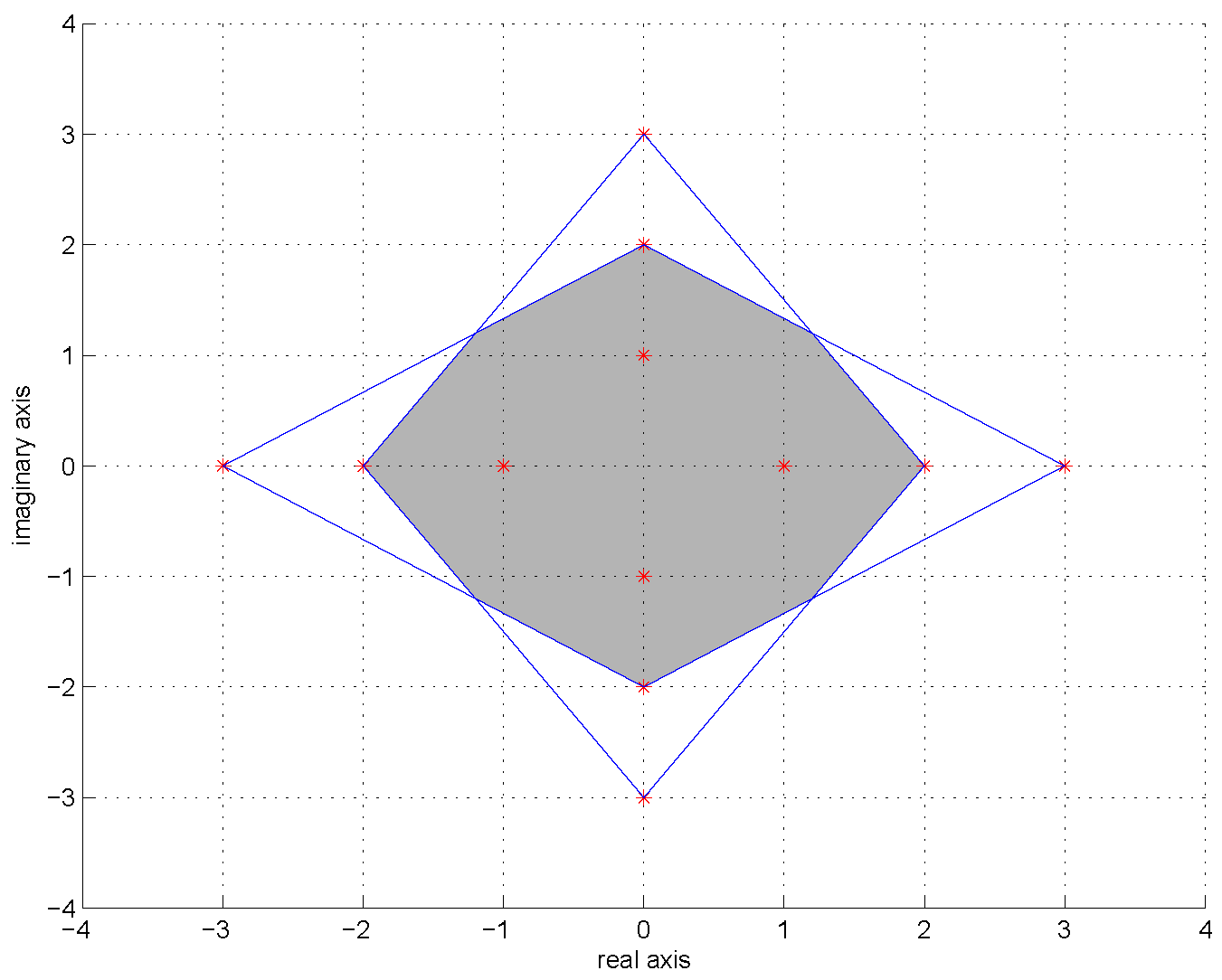, width =2in, height=2in} \qquad
\epsfig{file=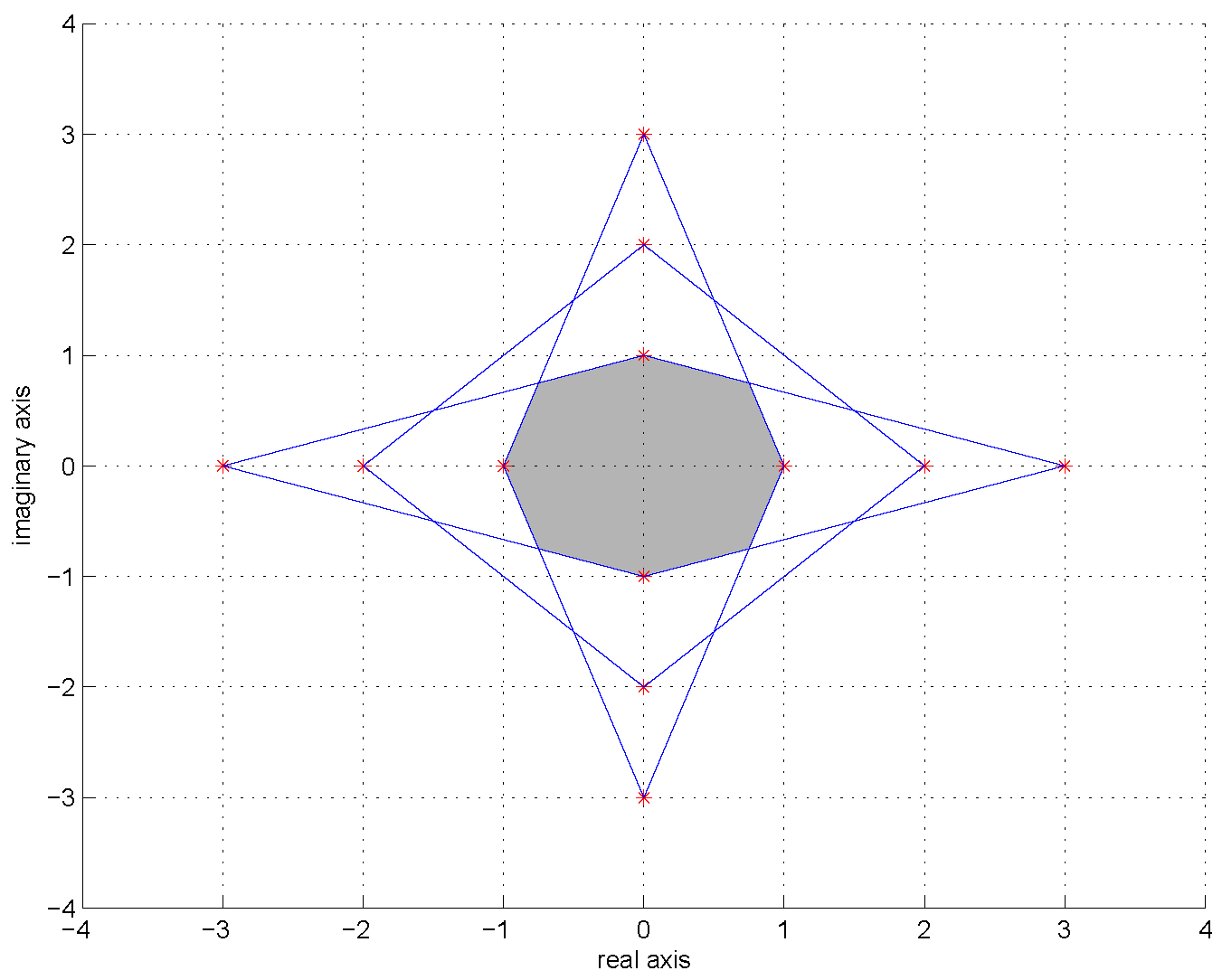, width =2in, height=2in}

$\Lambda_2(B)$\hspace{4.5cm}
$\Lambda_3(B)$

\end{center}

To deal with normal matrices $A \in M_n$ as in Example \ref{ex2} that
$\Lambda_1(A)$ is not an $n$-side convex polygon, we need to construct
some half spaces using the eigenvalues of the normal matrix $A$.
To do this, we introduce the following.
Given any two distinct complex numbers $a$ and $b$,
let $L(a,b)$ be the (directed) line passing through $a$ and $b$.
The closed half plane
$$H(a,b)
= \{ z \in \IC: \Im\left((\bar b-\bar a)(z-a)\right) \ge 0\}$$
is called the {\em left closed half plane} determined by $L(a,b)$.
For example, $H(0,i) = \{z\in \IC: \Re(z) \le 0\}$
and $H(i,0) = \{z\in \IC: \Re(z) \ge 0\}$.
Remark that in Example \ref{ex1}, the set $\cH_j$ is indeed the closed half plane $H(w^j,w^{j+k})$.
Note that $L(a,b) \ne L(b,a)$.
In our discussion, it is sometimes convenient to write
$$H(a,b) = \{z \in \IC: \Re( e^{-i\xi} z ) \le \Re( e^{-i\xi} a)  \}$$
with $\xi = \arg(b - a) - \pi/2$. Also, we use $H_0(a,b)$ to denote
the {\em left open half plane} determined by $L(a,b)$, i.e.,
$H_0(a,b) = H(a,b) \setminus L(a,b)$.

\smallskip
We have the following result showing that for a normal matrix $A \in M_n$
with $m$ distinct eigenvalues, $\Lambda_k(A)$ can be written as the
intersection of at most $\max\{m,4\}$ half spaces.
Even without any knowledge about the final shape of the set $\Lambda_k(A)$,
one can use $m(m-1)$ half spaces to generate $\Lambda_k(A)$.
Evidently, the construction is more efficient than the construction
using (\ref{eq1.1}) or (\ref{eq1.2}). Furthermore, we can conclude that
$\Lambda_k(A)$ is either an empty set, a singleton, a line segment,
or a non-degenerate polygon with at most $m$ vertices.

\smallskip
\begin{theorem} \label{thm2.3}
Let $A \in M_n$ be normal with distinct eigenvalues
$a_1, \dots, a_m$ that are not collinear.
Let $\cS$ be  the set of index pairs $(r,s)$ such that $H(a_r,a_s)$
contains at least $n-k+1$ eigenvalues (counting multiplicities) of $A$, and
\begin{eqnarray*}
\cS_0=\{&& (r,s)\in\cS: H_0(a_r,a_s)
\hbox{ contains at most } \cr
&&\hspace{3cm} n-k-1 \hbox{ eigenvalues (counting multiplicities) }\}.
\end{eqnarray*}
Then
\begin{eqnarray}\label{int}
\Lambda_k(A) = \bigcap_{(r,s) \in \cS} H(a_r,a_s)
=\bigcap_{(r,s) \in \cS_0} H(a_r,a_s).
\end{eqnarray}
Moreover, $\Lambda_k(A)$ can be written as intersection
of at most $\max\{m,4\}$ half planes $H(a_r,a_s)$, with $(r,s) \in \cS_0$.
\end{theorem}

\smallskip
\begin{proof}
In the first part of the proof,
we assume that $A\in M_n$ has $n$ eigenvalues $a_1,\dots,a_n$.
For notational simplicity, we write $H(a_r,a_s) = H(r,s)$,
$H_0(a_r,a_s) = H_0(r,s)$,
and $L(a_r,a_s) = L(r,s)$ for any two distinct eigenvalues
$a_r$ and $a_s$ of $A$.
For each $(r,s)\in \cS$, since $H(r,s)$ is convex and contains at least $n-k+1$
eigenvalues of $A$, by (\ref{eq1.2}), we have
$$
\Lambda_k(A) = \bigcap_{1 \le j_1<  \cdots < j_{n-k+1} \le n}
\conv\{a_{j_1}, \dots, a_{j_{n-k+1}}\}
\subseteq H(r,s).
$$
It follows that
\begin{equation}\label{eq2.2}
\Lambda_k(A)
\subseteq \bigcap_{(r,s)\in \cS} H(r,s).
\end{equation}

To prove the reverse inclusion of (\ref{eq2.2}), 
note that if $z$ is a point not in $\Lambda_k(A)$, then $z$ will lie outside
a convex polygon which equals the convex hull of  $n-k+1$ eigenvalues of $A$.  
So, it suffices to show that the convex hull $\cW$ of any $n-k+1$ eigenvalues of $A$ can be written
as an intersection of half planes, $\cW=\cap_{j=1}^\ell H(r_j,s_j)$ for some $(r_1,s_1),\dots, (r_\ell,s_\ell)\in\cS$. We consider the following three cases.

\smallskip\noindent
{\bf Case 1} Suppose $\cW$ is a singleton.
Then $\cW = \{a_r\}$ for some eigenvalue $a_r$
with multiplicity at least $n-k+1$.
Since the eigenvalues of $A$ are non-collinear,
there are eigenvalues $a_s$ and $a_t$ such that
$a_r$, $a_s$, and $a_t$ are not collinear. Then
$$\cW  = H(r,s) \cap H(s,r) \cap H(r,t) \cap H(t,r).$$

\smallskip\noindent
{\bf Case 2} Suppose $\cW$ is a non-degenerate line segment.
In this case, $\cW = \conv\{a_r,a_s\}$ for some eigenvalues $a_r$
and $a_s$ with $a_r \ne a_s$. Since the eigenvalues of $A$ are
non-collinear,  there is another eigenvalue $a_t$ such that
$a_r$, $a_s$, and $a_t$ are not collinear.
Without loss of generality,
we assume that $a_t \in H(r,s)$. Otherwise,
we interchange $a_r$ and $a_s$.
Then
$$\cW = H(r,s) \cap H(s,r) \cap H(s,t) \cap H(t,r).$$

\smallskip\noindent
{\bf Case 3} Suppose $\cW$ is a non-degenerate polygonal disk.
We may relabel the eigenvalues of $A$ and assume that
$\cW$ has vertices $a_1, \dots, a_q$ arranged in the
counterclockwise direction, where $q \ge 3$. For convenience of
notation, we will let $a_{q+1} = a_1$ and $H(q,q+1) = H(q,1)$.
Then $$\cW = \bigcap_{1\le t \le q }  H(t,t+1).$$
Thus, the first equality in (\ref{int}) is proved.
To prove the second equality in (\ref{int}),
we claim the following.

\smallskip
\noindent
{\bf Claim} For each $(r,s) \in \cS \setminus \cS_0$,
there exist two ordered pairs $(r_1,s_1)$ and $(r_2,s_2)$ in $\cS_0$
such that
$H(r_1,s_1) \cap H(r_2,s_2) \subseteq H_0(r,s)$.

\smallskip
Once the claim is proved,
all the half planes $H(r,s)$ with $(r,s) \in \cS \setminus \cS_0$
are not needed in the intersection $\bigcap_{(r,s) \in \cS} H(a_r,a_s)$
and hence the second equality in (\ref{int}) holds.

To prove the claim, suppose $(r,s)\in \cS\setminus \cS_0$.
Then $H_0(r,s)$ contains at least $n-k$ eigenvalues of $A$.
By a translation   followed by a rotation, we may assume that
$H(r,s) = \{z\in \IC: \Im z \ge 0\}$
and we can relabel the index of eigenvalues so that
for $1\le j\le n-1$,
either $\Im a_j > \Im a_{j+1}$
or $\Im a_j = \Im a_{j+1}$ with $\Re a_j \ge \Re a_{j+1}$. Let
$$\cU = \conv \{a_1,\dots, a_{n-k}\}
\quad\hbox{and}\quad
\cV = \conv \{a_{n-k+1}, \dots, a_n\}.$$
Then $\cU$ and $\cV$ are disjoint if $a_{n-k} \ne a_{n-k+1}$
or $\cU \cap \cV = \{a_{n-k}\}$ if $a_{n-k} = a_{n-k+1}$.
By the assumption, $\cU \subseteq H_0(r,s)$ and $\{a_r,a_s \} \subseteq \cV$.
Define the set
$$\cW = \conv \{a_i - a_j: 1\le i \le n-k < j \le n\}
= \{u-v: u\in \cU \hbox{ and } v \in \cV\},$$
which is a convex polygon.
Note that $\cW\subseteq\{z\in\IC:\Im(z)\ge 0\}$ since $\Im(a_i-a_j)\ge 0$ for all
$1\le i\le n-k <j\le n$.
By the facts that $\cU$ and $\cV$
can intersect at at most one point, and
the union $\cU\cup \cV$ cannot be contained in any line,
the set $\cW$ does not lie in any line that passes through the origin,
and the point $0$ can only be either an extreme point of $\cW$ or is not in $\cW$.
Under these conditions, one can find two extreme points
$w_1$ and $w_2$ in $\cW$ with $\Im (\bar w_1 w_2) \ne 0$ such that
\begin{eqnarray}\label{eq2.3}
\Im (\bar w_1 w) \ge 0 \ge \Im (\bar w_2 w) \quad\hbox{for all}\quad w\in \cW.
\end{eqnarray}
Since $w_1$ is an extreme point in $\cW$,
there are eigenvalues $a_{s_1} \in \cU$ and $a_{r_1} \in \cV$
such that $w_1 = a_{s_1} - a_{r_1}$.
Then (\ref{eq2.3}) gives
$$\Im (\bar a_{s_1} - \bar a_{r_1})(u - a_{r_1}) \ge 0
\quad\hbox{and}\quad
\Im (\bar a_{r_1} - \bar a_{s_1})(v - a_{s_1}) \ge 0$$
for all $u \in \cU$ and $v \in \cV$,
and thus, $\cU \subseteq H(r_1,s_1)$
and $\cV \subseteq H(s_1,r_1)$.
With the fact that $a_{r_1}$ and $a_{s_1}$ lie in the line $L(r_1,s_1)$,
the closed half plane $H(r_1,s_1)$ contains at least $n-k+1$ eigenvalues of $A$
while the open half plane $H_0(r_1,s_1)$ contains at most $n-k-1$ eigenvalues only.
Therefore, $(r_1,s_1) \in \cS_0$.
By a similar argument, one can show that there are eigenvalues $a_{r_2}\in \cU$ and $a_{s_2} \in \cV$
such that $w_2 = a_{r_2} - a_{s_2}$. Then (\ref{eq2.3}) yields
$$\Im (\bar a_{r_2} - \bar a_{s_2})(u - a_{s_2}) \le 0
\quad\hbox{and}\quad
\Im (\bar a_{s_2} - \bar a_{r_2})(v - a_{r_2}) \le 0$$
for all $u \in \cU$ and $v \in \cV$,
and thus, $\cU \subseteq H(r_2,s_2)$ and $\cV \subseteq H(s_2,r_2)$,
and one can conclude that $(r_2,s_2) \in S_0$.
Observe that the two lines $L(r_1,s_1)$ and $L(r_2,s_2)$ are not parallel
as $\Im (\bar a_{r_1} - \bar a_{s_1}) (a_{s_2} - a_{r_2} )
= \Im (\bar w_1 w_2)  \ne 0$.
Using the fact that the two distinct eigenvalues $a_r$ and $a_s$
are in $\cV$, which is contained in
the intersection $H(s_1,r_1) \cap H(s_2,r_2)$,
one can conclude that the intersection $H(r_1,s_1) \cap H(r_2,s_2)$
must lie in $H_0(r,s)$, the interior of $H(r,s)$.
Therefore, the claim holds.

Next, we turn to the last part of the Theorem.
It is trivial that if $\Lambda_k(A)$ is either an empty set, a singleton, or
a non-degenerate line segment, then only at most $4$ half planes are needed
in the construction of $\Lambda_k(A)$.

Suppose $A$ has $m$ distinct eigenvalues $a_1,\dots,a_m$ and
$\Lambda_k(A)$ is a non-degenerate polygon.
Let $\cT$ be a minimal subset of $\cS_0$ such that $\Lambda_k(A)=\cap_{(r,s)\in \cT}H(r,s)$.
Since $\cT$ is minimal, the half planes $H(r,s)$, $(r,s) \in \cT$, are all distinct.
We may further assume that for all $(r,s) \in \cT$, $\{ a_1,\dots, a_m\}\cap L(a_r,a_s)\subseteq \conv\{a_r,a_s\}$.
Since $\Lambda_k(A)$ is a non-degenerate polygon, for each $1\le t\le m$,
there exist at most two pairs $(r,s)\in \cT$ such that $t\in \{r,s\}$.
Therefore, $\cT$ contains at most $m$ ordered pairs.
\end{proof}

\smallskip
\begin{example} \rm
Let $A = \diag(0,0,1,1,i)$. Then
$$\begin{array}{rl}\Lambda_2(A)&=[0,1]=H(0,1)\cap H(1,0)\cap H(1,i)\cap H(i,0)\\&\\
\Lambda_3(A)&=\emptyset=H(1,0)\cap H(1,i)\cap H(i,0)\end{array}$$ and
the intersection of any $2$ half planes $H(a_r,a_s)$ is non-empty.
This example also shows that one cannot replace $\max\{m,4\}$ by $m$
in the conclusion in Theorem \ref{thm2.3}.
\end{example}

\smallskip
\begin{example} \rm
Let $A = \diag(1, -1, i, -i)$. Then
$$\Lambda_2(A)=\{0\}=H(1,-1)\cap
H(-1,1)\cap  H(i,-i)\cap H(-i,i)$$
and $\Lambda_2(A)$ cannot be written as an
intersection of less than 4
half planes $H(a_r,a_s)$.
\end{example}

\smallskip
\begin{corollary} \label{cor2.6}
Suppose $A \in M_n$ is normal such that $W(A)$ is an $n$-sided
polygon
containing the origin as its interior point. Let $v_1, \dots, v_n$ be the
vertices of $W(A)$ having arguments $0 \le \xi_1 < \cdots < \xi_n < 2\pi$. If
$k < n/2$, then  $\Lambda_k(A)$ is an $n$-sided convex
polygon  obtained by
joining $v_j$ and $v_{j+k}$, where $v_{j+k}=v_{j+k-n}$  if $j+k>n$.
\end{corollary}

\smallskip
By Theorem \ref{thm2.3}, it is easy to see that the boundary
of $\Lambda_k(A)$
are subsets of the union of line segments of the form
$\conv\{a_r,a_s\}$ such that $a_r$ and $a_s$ satisfy the $H(a_r,a_s)$ condition.
However, it is not easy to determine which part of the line segment
actually belong to $\Lambda_k(A)$ as shown in Examples \ref{ex1}, \ref{eg1},
and \ref{ex2}.
By Theorem \ref{thm2.3}, if the normal matrix $A \in M_n$
has $m$ distinct eigenvalues, we need no more than $\max\{m,4\}$
half planes $H(a_r,a_s)$ to generate $\Lambda_k(A)$.
{\it Can one determine these half planes effectively?}
We will answer this question by presenting an
algorithm in Section 5
based on the discussion in this section.

\section{Matrices with prescribed higher rank numerical ranges}

We study the following problem in this section.

\smallskip
\begin{problem} \label{prob.1} \rm
Let $k > 1$ be a positive integer,
and let $\cP$ be a $p$-sided   polygon in $\IC$.
Construct a normal matrix $A $ with smallest size (dimension) such that
$\Lambda_k(A) = \cP$.
\end{problem}

\smallskip
If $\cP$ degenerates to a line segment joining two points $a_1$ and $a_2$.
Then the smallest $n$ to get a normal matrix with $\Lambda_k(A)=\cP$
is $n=2k$, if $a_1$ and $a_2$ are distinct and $n=k$ if $a_1=a_2$.
So we focus on the case when the polygon $\cP$ is non-degenerate.

A natural approach to Problem \ref{prob.1} is to reverse
the construction of $\Lambda_k(A)$ in Example \ref{eg1}.
Suppose we have a non-degenerate $p$-sided polygon
$\cP$, with vertices, $v_1,\dots,v_p$. 

Without loss of
generality, we may assume that $0$ lies in the interior of
$\cP$ and the arguments of $v_j$ in $[0, 2\pi)$ are arranged in
ascending order. Our goal is to use the support line $L_j$
which passes through $v_j,v_{j+1}$ for $j = 1, \dots, p$, where
$v_{p+1} = v_1$,  to construct $A = \diag(a_1, \dots, a_p)$
such that $\Lambda_k(A) = \cP$. Note that if the desired values
$a_1, \dots, a_p$ exist and are arranged in counter-clockwise
direction, then (by proper numbering) the line $L_j$ will coincide
with the line passing through  $a_j$ and $a_{j+k}$,
where $a_{j+k}=a_{j+k-p}$ if $j+k>p$.
Consequently, $a_j$ will lie at the intersection of $L_j$ and
$L_{j-k}$, where $L_{j-k} = L_{j-k+p}$ if $j-k<0$.
Consequently, there exists $A = \diag(a_1, \dots, a_p)$ satisfying
$\Lambda_k(A) = \cP$ if the following hold.

\begin{itemize}
\item[{\rm (1)}]  $k < p/2$.

\item[{\rm (2)}] There exist $a_1, \dots, a_p \in \IC$ such that

\begin{itemize}
\item[{\rm  (2.a)}]  $L_j \cap L_{j-k} = \{a_j\}$ for $j = 1, \dots, p$,

\item[{\rm (2.b)}]  $a_1, \dots, a_p$ have arguments $\xi_1 < \dots < \xi_p$
in the interval $[\xi_1, \xi_1+2\pi)$ and 0 lie in the interior of
their convex hull.
\end{itemize}
\end{itemize}

\smallskip
Note that  by Theorem \ref{thm2.3},
$A$ has the smallest dimension among all
normal matrices $B$ such that $\Lambda_k(B) = \cP$.

Clearly, conditions (1) and (2a) are necessary in the above construction.
From the following example, one can see that
the above construction also fails when condition (2b)
is not satisfied.

%\smallskip
\begin{example} \rm
Let $\cP$ be the $5$-sided polygon with vertices
$\{v_1,\dots,v_5\} = \{2+i,1+2i,-1+3i,-1-i,3-i \}$, see the below figure.
\begin{center}
\epsfig{file=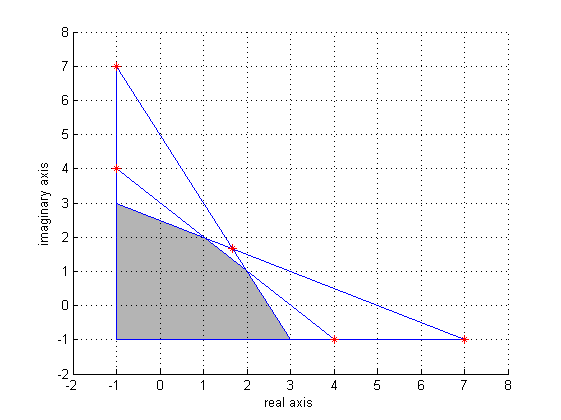, width =2in, height=2in} \\
The polygon $\cP$
\end{center}
Then, with $k=2$, we have  $$\{a_1,\dots,a_5\} = \{-1+4i,7-i,-1+7i,4-i, (5+5i)/3   \},$$
which does not satisfy the condition (2b). Clearly,
for $A = \diag(a_1,\dots,a_5)$, $\Lambda_2(A)$
lies in the convex hull of $\{a_1,\dots,a_5\}$,
which does not contain $\cP$.
\end{example}

\smallskip
Conditions (2a) and (2b) motivate the following definition.

\smallskip
\begin{definition} \label{3.1}
Let $\Omega = \{z\in \IC: |z| = 1\}$. A subset $\Pi = \{\alpha_1,\dots,\alpha_m\}$,
with distinct $\alpha_1,\dots,\alpha_m \in \Omega$,
is {\it $k$-regular} if every semi-circular arc of $\Omega$ without
endpoints contains at least $k$ elements in $\Pi$.
\end{definition}

\smallskip
Given distinct $\alpha_1,\ \alpha_2\in \Omega$, $\alpha_2/\alpha_1
=e^{i\theta}$
for a unique $0<\theta<2\pi$. Then
$[\alpha_1,\alpha_2]=\{e^{it}\alpha_1:0\le t\le \theta\} $ is the
closed arc on $\Omega$ from $\alpha_1$ to $\alpha_2$ in the
counterclockwise direction. Also define the open arc
$$(\alpha_1,\alpha_2)=[\alpha_1,\alpha_2 ]\setminus\{\alpha_1,\alpha_2\}.$$
The value $\theta$ is called the {\it length} of these intervals.
Suppose $1\le k\le n$ and $\Pi \subseteq \Omega$.
Then $\Pi$ is $k$-regular if for each
$\alpha\in \Pi$, $(\alpha,-\alpha)\cap \Pi $ contains at least $k$ elements.

Note that if
$\Pi=\{e^{i\xi_j}:1\le j\le n\}$ with distinct $\xi_1, \dots, \xi_n \in [0,2\pi)$,
then $\Pi$ is $k$-regular if and only if
for each $r = 1,\dots, n$, there are $1\le j_1 < \dots < j_k \le n$ such that
\begin{eqnarray}\label{eq3.1}
e^{i\xi_{j_1}},\dots,e^{i\xi_{j_k}} \in \left( e^{i\xi_r}, e^{i(\xi_r + \pi)} \right).
\end{eqnarray}
For this reason, a set
$\{\xi_1 , \dots , \xi_n\}\subseteq [0,2\pi)$ of $n$ distinct numbers is also
called $k$-regular if $\{e^{i\xi_j}:1\le j\le n\}$
is $k$-regular as defined in Definition \ref{3.1}.
For $\xi,\xi'\in [0,
2\pi)$, $[\xi,\xi']$ will denote the subset $\{t\in [0,
2\pi):e^{it}\in  [e^{i\xi},e^{i\xi'} ]\}$; the intervals $[\xi,\xi')$,
$(\xi,\xi']$ and $(\xi,\xi')$ will also be defined similarly.

In Example \ref{ex1}, a direct computation shows that for
$1\leq r, k\leq n$,
$$\xi_{r+k}-\xi_r = \cases{  
2k\pi/n  & \mbox{ if }  $r+k\leq n$, \cr
2k\pi/n-2\pi  & \mbox{ if } $r+k>n$.}
$$
Therefore, the set
$\{  \xi_1 ,\dots,  \xi_n \}$ is $k$-regular and
$\Lambda_k(A)$ is nonempty for $1\leq k< n/2$.
Otherwise, the set  $\{  \xi_1 ,\dots,  \xi_n \}$ is not
$k$-regular and $\Lambda_k(A)$ is either empty or a singleton.

In the following, we need an alternate formulation of
(\ref{eq1.1}). For any $d, \xi \in \IR$,
consider the closed half plane
\begin{equation} \label{chdxi}
\cH(d,\xi) = \{\mu\in \IC: \Re(e^{-i\xi} \mu ) \le d\},
\end{equation}
and its boundary, which is the
straight line
\begin{equation} \label{cldxi}
\cL(d,\xi) = \partial \cH(d,\xi) = \{\mu\in \IC: \Re(e^{-i\xi} \mu) = d\}.
\end{equation}
For $A\in M_n$, let Re$A=(A+A^*)/2$.
Then (\ref{eq1.1}) is equivalent to $$\Lambda_k(A) =
\bigcap_{\xi \in [0,2\pi)}
\cH(\lambda_k(\Re(e^{-i\xi}A)), \xi).$$
The following result is easy to verify.

\smallskip
\begin{proposition}
\label{prop3.a}
Let $A \in M_n$ and $\Lambda_k(A) = \cap_{j=1}^m \cH(d_j,\xi_j) \ne
\emptyset$,
where $\cH(d_j,\xi_j)$ is defined as in
$(\ref{chdxi})$
for some $d_1, \dots, d_m \in \IR$ and distinct $\xi_1, \dots, \xi_m\in   [0,2\pi)$.
\begin{itemize}
\item[{\rm (a)}] We have $0 \in \Lambda_k(A)$ if and only if
$d_1, \dots, d_m \ge 0$, and $0$ is an interior point of
$\Lambda_k(A)$ if and only if $d_1, \dots, d_m > 0$.
\item[{\rm (b)}] If $\mu = re^{i\xi}$ with $r > 0$ and $\xi \in \IR$,
then $\Lambda_k(\mu A) = \cap_{j=1}^m \cH(rd_j,\xi_j+\xi)$
and $\Lambda_k(A + \mu I ) = \cap_{j=1}^m \cH(d_j + r\cos(\xi - \xi_j),\xi_j)$.
\end{itemize}
\end{proposition}

\smallskip
In connection to Problem \ref{prob.1}, we have the following.

\smallskip
\begin{theorem}\label{equiv}
Suppose $\cP = \bigcap_{j=1}^p\, \cH(d_j,\xi_j)$
is a non-degenerate $p$-sided polygon,
where $\cH(d_j,\xi_j)$ is defined as in $(\ref{chdxi})$
with $d_1,\dots,d_p\in \IR$ and  distinct $\xi_1,\dots,\xi_p\in [0,2\pi)$.
Let $q$ be a nonnegative integer.
The following two statements are equivalent.
\begin{itemize}
\item[\rm (I)] There is a $(p+q)\times (p+q)$
normal matrix $A$ such that $\Lambda_k(A) = \cP$.
\item[\rm (II)] There are  distinct $\xi_{p+1},\dots,\xi_{p+q}\in [0,2\pi)$ such that
$\{\xi_1,\dots,\xi_{p+q}\}$ is $k$-regular.
\end{itemize}
\end{theorem}

\smallskip
Notice that a necessary condition for
the set $\bigcap_{j=1}^p\, \cH(d_j,\xi_j)$
to be a non-degenerate polygon is that
\begin{eqnarray}
\{e^{i\xi_1},\dots,e^{i\xi_p} \} \ \mbox{is $1$-regular.}
\end{eqnarray}
By Proposition \ref{prop3.a}, one may assume that $0$ lies in the interior
of $\cP$ in our proofs. However, it is equally convenient for us not
to impose this assumption so that we need not verify
$d_j > 0$ in $\cH(d_j,\xi_j)$ in our proofs.

To prove Theorem \ref{equiv}, we need some lemmas.

\smallskip
\begin{lemma}\label{lem3.5}
Given $A = \diag(a_1,\dots,a_n)$ and $1\le m < n$.
Suppose the eigenvalues $a_{m+1},\dots,a_n$ are in $\Lambda_k(A)$ but
not extreme points of $\Lambda_k(A)$. Then
$$\Lambda_k(\diag(a_1,\dots,a_m)) = \Lambda_k(A).$$
\end{lemma}

\smallskip
\begin{proof}
It suffices to show that if $a_n$ is in $\Lambda_k(A)$
but not an extreme point of $\Lambda_k(A)$, then
$\Lambda_k(\diag(a_1,\dots,a_{n-1})) = \Lambda_k(A)$.

Suppose $a_n$ satisfy the above assumption.
Clearly, $\Lambda_k(\diag(a_1,\dots,a_{n-1}))$
is a subset of $\Lambda_k(A)$.
On the other hand, for any $1\le j_1 < \dots < j_{n-k} \le n-1$,
$\Lambda_k(A)\subseteq\conv\{a_{j_1},\dots, a_{j_{n-k}}, a_n\}$.
Since $a_n$ is not an extreme point of $\Lambda_k(A)$,
it follows that $a_n$ lies in $\conv \{a_{j_1},\dots, a_{j_{n-k}},a_n\}$
but is not its extreme point. Therefore,
$$\conv \{a_{j_1},\dots, a_{j_{n-k}}\}
= \conv \{a_{j_1},\dots, a_{j_{n-k}},a_n\}.$$
Thus,
\begin{eqnarray*}
\Lambda_k(A)
&=& \bigcap \{ \conv \{a_{j_1},\dots, a_{j_{n-k+1}} \}:
1\le j_1 < \dots < j_{n-k} < j_{n-k+1}\leq n\} \cr
&\subseteq& \bigcap \{ \conv \{a_{j_1},\dots, a_{j_{n-k}},a_n \}:
1\le j_1 < \dots < j_{n-k} \le n-1\} \cr
&=& \bigcap \{ \conv \{a_{j_1},\dots, a_{j_{n-k}} \}:
1\le j_1 < \dots < j_{n-k} \le n-1\} \cr
&=& \Lambda_k(\diag(a_1,\dots,a_{n-1})).
\end{eqnarray*}
\vskip-.25in \hspace{55mm}
\end{proof}

\smallskip
The next lemma shows that if a convex polygon $\cP$ is the intersection
of half planes $\cH(d_j, \zeta_j)$ for $j = 1, \dots, m$,
such that the set $\{\zeta_1,\dots,\zeta_m\}$ is ``almost'' $k$-regular
(in the sense that $\{\zeta_1,\dots,\zeta_m\}$ is $k$-regular if
we count the multiplicity of each element in the set),
one may replace these half planes by $n$ other half planes
$\cH(\tilde d_j, \tilde \zeta_j)$
for $j = 1, \dots, n$, $n \ge m$, 
with  $\tilde \zeta_i \ne \tilde \zeta_j$
for all $i\ne j$, such that 
 $\{\tilde \xi_1,\dots,\tilde \xi_n\}$  
is $k$-regular and the boundary $\cL(\tilde d_j, \tilde \zeta_j)$ of
$\cH(\tilde d_j, \tilde \zeta_j)$ touches the polygon $\cP$
for each  $j = 1, \dots, n$.

\smallskip
\begin{lemma}\label{lem3.6}
Suppose $\cP = \cap_{j=1}^m \cH(d_j, \zeta_j)$ such
that $0 \le \zeta_1 \le \dots \le \zeta_m < 2\pi$ and for each $r = 1,\dots,m$,
there are $1 \le j_1 < \cdots < j_k \le m$
such that $\zeta_{j_1}, \dots, \zeta_{j_k} \in (\zeta_r, \zeta_r + \pi)$.
For every $n \ge m$, there exist
$\tilde d_1, \dots, \tilde d_n\in\IR$ and
distinct $\tilde \zeta_1,\dots,\tilde \zeta_n \in [0,2\pi)$
with
$\{\tilde \zeta_1, \dots, \tilde \zeta_n\}$
being $k$-regular
such that
$\cP = \cap_{j=1}^n \cH(\tilde d_j, \tilde \zeta_j)$ and
$\cP \cap \cL(\tilde d_j, \tilde \zeta_j) \ne \emptyset$
for each $j = 1, \dots, n$.
\end{lemma}

\smallskip
\begin{proof}
Set $\tilde \zeta_1 = \zeta_1$, and
for $s\in\{2, \ldots, m\}$,  let $\tilde \zeta_{s}=\zeta_{s}$
if $\zeta_{s-1}<\zeta_s$.  For the remaining values, we have
$\zeta_{s-1}=\zeta_s$  and we can set
$$
\tilde \zeta_{s-t_1}=\zeta_{s-t_1+1}=\cdots=\zeta_s =\cdots=
\zeta_{s+t_2}<\tilde\zeta_{s+t_2+1}$$ for some $t_1\ge 1$ and
$t_2\ge 0$. Let $\ell=\min\{j:\tilde\zeta_j>\tilde\zeta_{s-t_1}+\pi\}$, then we can
replace $\zeta_{s+j}$ by $\tilde
\zeta_{s+j}=\zeta_{s+j}+\epsilon_j$ for sufficient small
$\epsilon_j>0$ for $j=-t_1+1, -t_1+2, \ldots, 0, \ldots, t_2$ such
that $$
\tilde \zeta_{s-t_1}<\tilde\zeta_{s-t_1+1}<\cdots<\tilde\zeta_s <\cdots<
\tilde\zeta_{s+t_2}<\min\{\tilde\zeta_\ell-\pi, \tilde\zeta_{s+t_2+1}\}.$$
After this modification, $\tilde\zeta_1, \ldots,
\tilde\zeta_m$ are distinct and $\{\tilde\zeta_1, \ldots, \tilde\zeta_m\}$
is $k$-regular. If $n > m$, pick distinct $\tilde \zeta_{m+1},\dots,\tilde \zeta_n
\in [0,2\pi) \setminus \{\tilde \zeta_1,\dots,\tilde \zeta_m\}$.
Then $\{\tilde \zeta_1,\dots,\tilde \zeta_n\}$ also forms a $k$-regular set.
Finally, let $\tilde d_j = \max_{\mu\in \cP}\ \Re
\left( e^{-i\tilde \zeta_j} \mu \right)$ for $j =1,\dots,n$. Clearly,
we have $\cP \cap \cL(\tilde d_j, \tilde \zeta_j) \ne \emptyset$ and
$\cP\subseteq\cH(\tilde d_j, \tilde \zeta_j)$ for all $j$. By construction,
$\{\zeta_1, \ldots, \zeta_m\}\subseteq  \{\tilde\zeta_1, \ldots,
\tilde\zeta_n\}$  and $\cP = \cap_{j=1}^m \cH(d_j, \zeta_j)=
\cap_{j=1}^n \cH(\tilde d_j, \tilde \zeta_j)$.

\smallskip
We can now present the proof of Theorem \ref{equiv}.

\smallskip
\bf Proof of Theorem \ref{equiv}. \rm
Let $\cP = \bigcap_{j=1}^p\, \cH(d_j,\xi_j)$
be a non-degenerate $p$-sided polygon, where
$d_1,\dots,d_p\in \IR$ and $\xi_1,\dots,\xi_p\in [0,2\pi)$.

\smallskip
Suppose (I) holds. We may assume that
$A = \diag(a_1,\dots,a_{p+q})$ and $\Lambda_k(A) = \cP$.
By Lemma \ref{lem3.5}, one can remove the eigenvalues of $A$
in $\Lambda_k(A)$ that are not extreme points of $\Lambda_k(A)$
to get $\tilde A \in M_n$ for some positive integer $n \le p+q$
so that $\Lambda_k(A) = \Lambda_k(\tilde A)$.
We have the following.

\smallskip
\noindent
{\bf Claim}
There are $f_1,\dots,f_n \in \IR$ and   $\zeta_1, \dots, \zeta_n \in [0,2\pi)$
such that $\Lambda_k(\tilde A) = \linebreak
\cap_{j=1}^n \cH(f_j, \zeta_j)$. Furthermore, for each $r = 1,\dots,n$, there exist $1 \le j_1 < \dots < j_k \le n$ such that
$\zeta_{j_1},\dots,\zeta_{j_k} \in (\zeta_r, \zeta_r + \pi)$.

\smallskip
Once the claim holds, Lemma \ref{lem3.6} will ensure that
$\Lambda_k(\tilde A) = \cap_{j=1}^{p+q} \cH(\tilde d_j, \tilde
\xi_j)$ for some $\tilde d_1,\dots,\tilde d_{p+q} \in \IR$ and a
$k$-regular set  $\{\tilde  \xi_1, \dots, \tilde \xi_{p+q}\}$,   with
$$ \bigcap_{j=1}^p\, \cH(d_j,\xi_j)
= \cP = \Lambda_k(A) = \Lambda_k(\tilde A) = \cap_{j=1}^{p+q} \cH(\tilde d_j, \tilde
\xi_j)\,.$$
Then
$\{\xi_1,\dots,\xi_p\} \subseteq \{\tilde \xi_1,\dots,\tilde \xi_{p+q}\}$. Thus,
we can take $\xi_{p+1},\dots,\xi_{p+q}\in [0,2\pi)$ so that
$\{\xi_1,\dots,\xi_{p+q}\}
= \{\tilde \xi_1,\dots,\tilde \xi_{p+q}\}$. Therefore, (II) holds.

For notational convenience, we assume that $A = \tilde A$ 
in the claim so that every eigenvalue of $A$ is
either an extreme point of $\Lambda_k(A)$ or does not lie in
$\Lambda_k(A)$.

We first construct $\zeta_1,\dots,\zeta_n \in [0,2\pi)$ and $f_1,\dots,f_n\in \IR$.
For $r = 1,\dots, n$,
let $\Gamma_r$ be the set containing all $\xi \in [0,2\pi)$
such that the closed half plane $\cH( \Re(e^{-i\xi} a_r), \xi)$ contains at least $n-k+1$ eigenvalues of $A$.
As $a_r$ is either an extreme point of $\Lambda_k(A)$ or not in
$\Lambda_k(A)$, there is a $\zeta\in [0,2\pi)$
such that $\Re(e^{-i\zeta} a_r) \ge \lambda_k(\Re(e^{-i\zeta} A))$
and hence the half plane $\cH( \Re(e^{-i\zeta} a_r), \zeta)$ contains at least $n-k+1$ eigenvalues.
Then $\Gamma_r$ is always nonempty.
Furthermore, by the definition of $\Gamma_r$,
the set $\Gamma_r$ is an union of closed arcs of $\Omega$.
Clearly,
$$\cP = \Lambda_k(A) \subseteq
\bigcap_{\xi \in \Gamma_r} \cH( \Re(e^{-i\xi} a_r), \xi).$$
Also the above intersection, which containing $\cP$, is a non-degenerate conical region.
Then $\Gamma_r$ is contained in some open semi-circular arc of $\Omega$; otherwise,
the above intersection of half planes is equal to the singleton $\{a_r\}$.
As $\Gamma_r$ is a union of closed arcs in some open semi-circular arc of $\Omega$,
there exists a unique $\zeta_r \in \Gamma_r$ such that
\begin{eqnarray}\label{eq3.7}
\Gamma_r \subseteq (\zeta_r - \pi,\zeta_r].
\end{eqnarray}

Let $f_r = \Re(e^{-i\zeta_r} a_r)$ for $1\le r \le n$.
We show that $\Lambda_k(A) =  \bigcap_{j = 1}^n \cH(f_j, \zeta_j)$.
Suppose $\cT$ is a minimal subset of $\cS_0$ such that $\Lambda_k(A) = \bigcap_{(r,s)\in \cT} H(a_r,a_s)$.
We may further assume that for all $(r,s)\in \cT$,
$\{a_1,\dots,a_m\} \cap L(a_r,a_s) \subseteq \conv \{a_r,a_s\}$.
For each $(r,s) \in \cT$,
write $H(a_r,a_s) = \cH( \Re(e^{-i\zeta} a_r), \zeta)$
with $\zeta = \arg(a_s-a_r) - \pi/2$. Then $\zeta \in \Gamma_r$.
We claim that $\zeta = \zeta_r$.
Suppose not. By the above assumption on $a_r$,
one can see that for a sufficiently small $\epsilon > 0$,
the half plane $\cH( \Re(e^{-i\hat \zeta} a_r), \hat \zeta)$ with $\hat \zeta = \zeta - \epsilon$
will contain all eigenvalues of $A$ that are in $H(a_r,a_s)$,
i.e., $\hat \zeta \in \Gamma_r$. With (\ref{eq3.7}), we have
$$\Lambda_k(A) \subseteq
\cH( \Re(e^{-i\zeta_r} a_r), \zeta_r) \cap
\cH( \Re(e^{-i\hat \zeta} a_r), \hat \zeta) \subseteq
H_0(a_r,a_s)\cup \{a_r\}.$$
So  $L(a_r,a_s)\cap\Lambda_k(A)$ contains at most one point.
But this contradicts the fact that $(r,s)$ is an element in
the minimal subset $\cT$.
Therefore, $\zeta =\zeta_r$.
Then for each $(r,s)\in \cT$, $H(a_r,a_s) = \cH(f_r,\zeta_r)$ and so
$$\bigcap_{j = 1}^n \cH(f_j,\zeta_j)
\subseteq \bigcap_{(r,s)\in \cT} H(a_r,a_s)
= \Lambda_k(A) \subseteq \bigcap_{j = 1}^n \cH(f_j, \zeta_j).$$
Thus, $\Lambda_k(A) =  \bigcap_{j = 1}^n \cH(f_j, \zeta_j)$
and the first part of the claim holds.

To prove the second part of the claim, without loss of generality, we may assume that
$a_r = 0$ and $\zeta_r = 0$. Then $f_r = 0$ and
$$\cH(f_r,\zeta_r) = H = \{z\in \IC: \Re(z) \le 0\}.$$
Thus, the closed left half plane contains at least $n-k+1$ eigenvalues of $A$.
Suppose that the closed right half plane $-H$ contains eigenvalues
$a_{j_1},\dots,a_{j_{h}}$ of $A$ with
$\zeta_{j_t} \ne 0$ for $t = 1,\dots, g$, and
$\zeta_{j_t} = 0$ for $t = g+1,\dots,h$
for some $g \le h$.
Fix a sufficiently small $\epsilon > 0$.
We choose $g+1\le \ell\le h$ so that
$$\Re( e^{-i\epsilon} a_{j_\ell})
= \max_{g+1\le t \le h} \Re( e^{-i\epsilon} a_{j_t}).$$
Then $\{a_{j_{g+1}},\dots, a_{j_h}\} \subseteq
\cH(\Re( e^{-i\epsilon} a_{j_\ell}), \epsilon)$.
On the other hand, this closed half plane
$\cH(\Re( e^{-i\epsilon} a_{j_\ell}), \epsilon)$
also contains all eigenvalues of $A$ that are in the left open half plane.
Thus, this closed half plane $\cH(\Re( e^{-i\epsilon} a_{j_\ell}), \epsilon)$
has at least $n-g$ eigenvalues of $A$.
On the other hand by (\ref{eq3.7}),
$\epsilon \notin\Gamma_{j_\ell}$ and so
$\cH(\Re( e^{-i\epsilon} a_{j_\ell}), \epsilon)$ can have at
most $n-k$ eigenvalues. Thus, we have $g \ge k$.

Now for each $t=1,\dots,k$, let $\hat d_t =\Re(a_{j_t})$, then
$\hat d_t \ge 0$ and $H \subseteq \cH(\hat d_t,0)$.
Thus, the closed half plane $\cH(\hat d_t,0)$ contains at least $n-k+1$ eigenvalues of $A$,
i.e., $\zeta_r = 0 \in \Gamma_{j_t}$. Recall that $\zeta_{j_t} \ne 0$.
By (\ref{eq3.7}), one see that
\begin{eqnarray*}
\zeta_r \in (\zeta_{j_t} - \pi, \zeta_{j_t})
\quad \hbox{ for } t = 1,\dots,k.
\end{eqnarray*}
Equivalently, $\zeta_{j_1},\dots,\zeta_{j_k} \in (\zeta_r, \zeta_r+\pi)$.
Thus, our {\bf claim} is proved, and (II) holds.

\smallskip
Suppose now (II) holds, namely, there are distinct $\xi_{p+1},\dots,\xi_{p+q}$ such that \linebreak
$\{\xi_1,\dots,\xi_{p+q}\}$ is $k$-regular. For $j=p+1,\dots,p+q$, define
$$d_j = \max_{\mu\in \cP}\ \Re ( e^{-i\xi_j} \mu  ).$$
Then $\cP \subseteq \cH(d_j,\xi_j)$ and so
$$\cP = \bigcap_{j=1}^p\, \cH(d_j,\xi_j)
=  \bigcap_{j=1}^n\, \cH(d_j,\xi_j)$$
with $n = p+q$.
By Lemma \ref{lem3.6}, we may assume that
$\cP\cap \cL(d_j,\xi_j)\neq\emptyset$ for all $j = 1, \dots, n$,
and $0\le\xi_1< \cdots < \xi_n<2\pi$ such that
condition (\ref{eq3.1}) holds.
For each $r=1,\dots,n$, let
$$a_r = \frac{i}{\sin(\xi_{r+k} - \xi_r)}
\left(  e^{i\xi_r}d_{r+k}  -  e^{i\xi_{r+k}} d_r\right)$$
and $A = \diag(a_1,\dots,a_n)$. Then
$$\Re(e^{-i\xi_r} a_r) = d_r
\quad\hbox{and}\quad
\Re(e^{-i\xi_{r+k}} a_r) = d_{r+k}.$$
Note that $a_r \in \cL(d_r,\xi_r) \cap \cL(d_{r+k}, \xi_{r+k})$ is
the vertex of the conical region $\cH(d_r,\xi_r)\cap\cH(d_{r+k},
\xi_{r+k})$, which contains $\cP$.  Therefore, $$
\Re(e^{-i\xi_r} (a_r-\mu)) \ge 0\quad
\mbox{and}\quad
\Re(e^{-i\xi_{r+k}}( a_r-\mu))\ge 0,$$
for all $\mu\in\cP$. Since $\xi_{r+k} \in (\xi_r,\xi_r+\pi)$, we have
\begin{eqnarray}\label{eq3.5}
\Re(e^{-i\xi} a_r) \ge \max_{\mu \in \cP}\Re(e^{-i\xi} \mu)
\quad \hbox{for all }\xi \in [\xi_r,\xi_{r+k}].
\end{eqnarray}
Let $\mu_j\in
\cL(d_j,\xi_j)\cap\cP$ for $j=r,r+k$.
As $\xi_{r+k} \in (\xi_r,\xi_r+\pi)$, we have
$\mu_r=a_r-ie^{i\xi_r}b_r$ and $\mu_{r+k}=a_r+ie^{i\xi_{r+k}}c_r$
for some $b_r, c_r\ge 0$. Note that
$$\Re(e^{-i\xi}(\mu_r-a_r))=b_r\sin(\xi_r-\xi)\ge 0 \quad
\mbox{for all} \ \xi\in[\xi_r-\pi, \xi_r],$$
and $$\Re(e^{-i\xi}(\mu_{r+k}-a_r))=c_r\sin(\xi - \xi_{r+k})\ge 0 \quad
\mbox{for all} \ \xi\in[\xi_{r+k}, \xi_{r+k}+\pi].$$
Since $\{\xi_1, \ldots, \xi_n\}$ is $k$-regular,
it is easily seen that
$$[0,2\pi) \setminus [\xi_r,
\xi_{r+k}]=[\xi_r-\pi, \xi_r)\cup(\xi_{r+k}, \xi_{r+k}+\pi].$$
Therefore, for  $\xi \in [0,2\pi) \setminus [\xi_r,
\xi_{r+k}]$, we have
$$\max\{\Re(e^{-i\xi}(\mu_r-a_r)),\Re(e^{-i\xi}(\mu_{r+k}-a_r))\}\ge 0.$$
Moreover, we have
\begin{equation}\label{eq3.6}
 \max\{\Re(e^{-i\xi } \mu_r), \Re(e^{-i\xi }
\mu_{r+k})\}\ge \Re(e^{-i\xi} a_r).
\end{equation}

Let $\xi \in [0,2\pi)$. Then $\xi \in [\xi_s,\xi_{s+1})$ for some
$s \in\{1,\dots,n\}$. It follows that $\xi \in [\xi_r, \xi_{r+k}]$
for $r = s-k+1,\dots,s$, and $\xi \in [0,2\pi) \setminus [\xi_r,
\xi_{r+k}]$ for other $r$. By (\ref{eq3.5}) and (\ref{eq3.6}),
$$\min_{r \in \{s-k+1,\dots,s\}} \Re(e^{-i\xi} a_r)
\ge \max_{\mu \in \cP} \Re(e^{-i\xi} \mu)
\geq \max_{r \notin \{s-k+1,\dots,s\}} \Re(e^{-i\xi} a_r).$$
Thus, $\lambda_k( \Re(e^{-i\xi} A))
= \min_{r \in\{s-k+1,\dots,s\}} \Re(e^{-i\xi} a_r)$ and so
$$\cP \subseteq \cH\left(\lambda_k(\Re(e^{-i\xi} A) ),\xi \right).$$
Hence, $\cP \subseteq \Lambda_k(A)$.
Furthermore, if $\xi = \xi_s$, then $\Re(e^{-i\xi_s} a_s) = d_s$. Thus
$$\lambda_k(\Re(e^{-i\xi_s} A)) = \min_{r \in \{s-k+1,\dots,s\}}
\Re(e^{-i\xi_s} a_r) \le d_s.$$
It follows that
\begin{eqnarray*}
\Lambda_k(A)
&=& \bigcap_{\xi\in [0,2\pi)} \cH
\left( \lambda_k(\Re(e^{-i\xi} A)), \xi\right)   \\
&\subseteq& \bigcap_{1\le s \le n} \cH
\left( \lambda_k(\Re(e^{-i\xi_s} A)), \xi_s\right)
\subseteq \bigcap_{1\le s \le n} \cH \left(d_s, \xi_s \right)
= \cP.
\end{eqnarray*}
Thus, $\cP = \Lambda_k(A)$.
\end{proof}

\smallskip
By  Theorem \ref{equiv}, Problem \ref{prob.1} is equivalent
to the following combinatorial problem, whose solution will
be given in the next section.

\medskip
\begin{problem}\label{prob.2}
Suppose $\{\xi_1,\dots,\xi_p\} \subseteq [0,2\pi)$ is $1$-regular. For $k>1$,
determine the smallest nonnegative integer $q$ so that
$\{\xi_1,\dots,\xi_{p+q}\}$ is $k$-regular
for some distinct $\xi_{p+1},\dots,\xi_{p+q}\in [0,2\pi)$.
\end{problem}

\section{Solutions for Problems \ref{prob.1} and \ref{prob.2}}

In this section, we give the solutions for Problems
\ref{prob.1} and \ref{prob.2}.
Given a non-empty set
$\Pi = \{\xi_1,\dots,\xi_p\} \subseteq \Omega$,
Problem \ref{prob.2} is equivalent to the study
of smallest nonnegative integer $q$ so that $\{\xi_1,\dots,\xi_{p+q}\}$
is $k$-regular for some distinct
$\xi_{p+1},\dots,\xi_{p+q} \in \Omega$.
We have the following.

\medskip
\begin{theorem}\label{main6}
Let $k>1$ be a positive integer and $\Pi$ be a $p$ element subset
of $\Omega$, including $s$ pairs of antipodal points:
$\{\beta_1, -\beta_1\}, \dots, \{\beta_s,-\beta_s\}$,
where $p \ge 3$ and $s \ge 0$.
Suppose $\Pi$ is $1$-regular but not $k$-regular
and $q$ is the minimum number of
points in $\Omega$ one can add to $\Pi$ to form a $k$-regular set.

\noindent
{\rm (a)} If $k \ge p-s$, then
\begin{eqnarray}\label{Q1}
q =
\cases{
2k+1-p & \hbox{if } $s = 0$, \cr
2k+2-p & \hbox{if } $s > 0$.}
\end{eqnarray}

\noindent
{\rm (b)} If $k < p-s$, then $q$ is the smallest nonnegative integer $t$
such that one can remove $t$ non-antipodal points from $\Pi$  to get a $(k-t)$-regular set. More precisely, 
\begin{eqnarray} \label{Q2}
q  =  &&\min \{ t\in \IN:
\Pi \setminus \{ \beta_1, -\beta_1, \dots, \beta_s,-\beta_s \}
\mbox{ has a {\it t}-element }  \cr
&&\hspace{3cm} 
\mbox{ subset $T$ such that }
\Pi \setminus T \mbox{ is }
(k-t)\mbox{-regular}\, \}.
\end{eqnarray}
Consequently,
\begin{equation} \label{optimal1}
q \le \min\{2k+2-p,k-1\}.
\end{equation}
The inequality in (\ref{optimal1})
becomes equality if
$\Pi = \{1, i, -1,\alpha_4,\dots,\alpha_p\}$
where $\alpha_4,\dots,\alpha_p$ lie in the open lower half plane.
\end{theorem}

\smallskip
Several remarks concerning Theorem \ref{main6} are in order.
If condition (a) in the theorem holds, then the value $q$
can be determined immediately.  However, it is important to
consider two cases depending on whether $\Pi$ has pairs of antipodal
points as illustrated by the following.

\medskip
\begin{example} \rm
Suppose $S_1 = \{1, w, w^2, w^3\}$ with $w = e^{2i\pi/5}$.
Then $\alpha \ne -\beta$ for any two elements $\alpha, \beta \in S_1$
and adding  $w^4$ to $S_1$ results in a $2$-regular set.
Suppose $S_2 = \{1, -1, i, -i\}$. Then we need to add at least two
points, say, $z, -z \in \Omega \setminus S_2$, to get a $2$-regular set.
\end{example}

\smallskip
Suppose condition (b) in the theorem
holds.
We can determine the value $q$ by taking $t$
non-antipodal elements away from $\Pi$
at a time and check whether the resulting set is $(k-t)$-regular.
The value $q$ can then be determined in no more than $\sum_{i=0}^{p-2s}{\tiny\(\begin{array}{c}p-2s\\ i\end{array}\) } = 2^{p-2s}$ steps.
The success of reducing Problem \ref{prob.2} to a problem which is solvable
in finite steps depends on Lemma \ref{prop2} and Proposition \ref{prop3}.

It would be nice to have a simple formula for $q$ in terms of $p,k,s$
in case (b) of the theorem.
However, the following example show that the value
$q$ depends not only  on the values
$p$,$k$,$s$, but also   on the relative positions
of the points in $\Pi$.

\medskip
\begin{example} \rm
Let $S_1 = \{1,w,w^2,w^3,w^4,w^5\}$ with $w = e^{2\pi i/7}$
and $S_2 = \{z^2, z^3, z^7, z^8, z^{12}, z^{13}\}$
with $z = e^{2\pi i/15}$.
Notice that both $S_1$ and $S_2$ contain $6$ elements
and have no antipodal pairs.
Furthermore, both of them are $2$-regular but not $3$-regular.
Clearly, adding $w^6$ to $S_1$ results a $3$-regular set.
However, as each of the open arcs
$(z^3, -z^3)$, $(z^8,-z^8)$ and $(z^{13},-z^{13})$
contains only two elements of $S_2$ while
the intersection of this three open arcs is empty,
at least two elements has to be added
to $S_2$ to form a $3$-regular set.
\end{example}

\smallskip
Note that our proofs are constructive;
see Lemma \ref{prop2} and Propositions \ref{ass2} and \ref{prop3}.
One can actually construct a subset
$\Pi' \subseteq \Omega$ with $q$ elements so that $\Pi \cup \Pi'$ is
$k$-regular.

By Theorem \ref{main6}, we can
answer Problems \ref{prob.1} and \ref{prob.2}, and obtain some
additional information on their solutions. We will continue to use
the notation $\cH(d,\xi)$ defined in (\ref{chdxi}) in the following.

\medskip
\begin{theorem}\label{main7}
For Problem \ref{prob.1}, if a $p$-sided polygon $\cP$
is expressed  as $\cP = \cap_{j=1}^p \cH(d_j,\xi_j)$ for some
$d_1, \dots, d_p \in \IR$ and $\xi_1,\dots,\xi_p \in [0,2\pi)$,
then the minimum dimension $n$ for the existence of a normal matrix
$A \in M_{n}$ such that $\Lambda_k(A) = \cP$
is equal to $p+q$, where $q$ is determined in Theorem \ref{main6}.
Moreover,
\begin{equation} \label{optimal2}
n \le \max\{2k+2,p+k-1\}.
\end{equation}
The inequality in (\ref{optimal2})
becomes equality if
$(\xi_1, \xi_2, \xi_3) = (0, \pi/2, \pi)$
and $\xi_4,\dots,\xi_p$ lie in $(\pi,2\pi)$.
\end{theorem}

\smallskip
We break down the proofs of Theorems \ref{main6} and \ref{main7}
in several propositions.
We first give a lower bound for the number of elements in a $k$-regular set.

\medskip
\begin{proposition}\label{lower}
Suppose $S = \{\alpha_1,\dots,\alpha_n\} \subseteq \Omega$ is $k$-regular.
Then $n \ge 2k+1$. Furthermore, if
$S$ contains a pair of antipodal points $\{\alpha,-\alpha\}$,
then $n \ge 2k+2$.
\end{proposition}

\smallskip
\begin{proof}
For any $r \in \{1,\dots, n\}$, each of the open arcs
$(\alpha_r,-\alpha_r)$ and $(-\alpha_r,\alpha_r)$
contains $k$ elements of  $S$. Thus, $n \ge 2k+1$.
For the last statement, if we take $\alpha_r = \alpha$,
then together with $\alpha$ and $-\alpha$, we see that
$n \ge 2k+2$. The proof of the assertion is complete.
\end{proof}

\smallskip
As shown in Proposition \ref{lower}, the existence of
a pair of antipodal points $\{\alpha, -\alpha\}$
has implication on the size of a $k$-regular set $\Pi$.
The next result together with Proposition \ref{lower} show that
the lower bound in  (\ref{Q1}) is best possible.

\medskip
\begin{proposition}\label{ass2}
Let $k > 1$ and $\Pi$
is a $p$ element subset of $\Omega$ containing $s$ pairs of antipodal points,
where $p \ge 3$ and $s \ge 0$.
If $\Pi$ is $1$-regular but not $k$-regular and $k \ge p-s$,
then one can extend $\Pi$ to a $k$-regular set by adding $2k+1-p$ or $2k+2-p$ elements,
depending whether $s$ is zero.
\end{proposition}

\smallskip
\begin{proof}
Assume $k \ge p-s$.
Suppose first that
$s > 0$.
Let $\Pi''$ be a set containing $(k-p+s+1)$ pairs
of antipodal points such
that $\Pi'' \cap \Pi$ is empty. Take
$$\Pi' = \Pi'' \cup -(\Pi \setminus \{\beta_1,-\beta_1,\dots,\beta_s,-\beta_s\}).$$
Then $\Pi'$ contains $(2k+2-p)$ elements.
Furthermore, the set $\Pi\cup\Pi'$ contains exactly $k+1$ pairs
of antipodal points hence it is $k$-regular.
Thus, the result follows if $s > 0$.

Next, suppose $s=0$.
Without loss of generality, we may assume that
$1\in \Pi$. Hence,  $-1\in \Pi'$.
We now modify $\Pi'$.
We first delete the point $-1$ in $\Pi'$.
Then for all other points $\alpha \in \Pi'$,
we replace $\alpha$ by $e^{i\xi} \alpha$ if $\alpha$ lies in the upper open half plane
$P = \{z \in \IC: \Im(z) > 0\}$,
and by $e^{-i\xi} \alpha$ if $\alpha$ lies in the lower open half plane $-P$,
with sufficiently small   $\xi > 0$.
Then we see that for every $\alpha \in \Pi\cup \Pi'$,
$\alpha P$ still contains exactly $k$ elements.
Thus, $\Pi\cup \Pi'$ is $k$-regular. Furthermore,
the modified set $\Pi'$ has one fewer point, i.e.,
$\Pi'$ has only $2k+1-p$ elements.
The proof of is complete.
\end{proof}

\smallskip
A referee pointed out that each $(k-1)$-regular set can
be enlarged to a $k$-regular set by adding in not more than $2$ extra elements.
The following result shows that sometimes $2$ may not be the minimum number needed.

\medskip
\begin{lemma}\label{prop2}
Let $k > 1$ and $\Pi$ be a subset of $\Omega$
containing at least one non-antipodal point.
The following are equivalent.
\begin{enumerate}
\item[\rm (a)] One can add a point $\beta\notin \Pi$ so that
$\Pi \cup \{\beta\}$ is $k$-regular.
\item[\rm (b)] One can delete a non-antipodal point $\gamma \in \Pi$ so that
$\Pi \setminus \{\gamma\}$ is $(k-1)$-regular.
\end{enumerate}
Here, an element $\alpha\in \Pi$ is called a non-antipodal point of $\Pi$ if $-\alpha \notin \Pi$.
\end{lemma}

\smallskip
\begin{proof}
Suppose first that (b) holds.
Let $P = \{z\in \IC: \Im(z) > 0\}$.
Without loss of generality, we may assume that $\gamma = 1$
is a non-antipodal point in $\Pi$. Suppose $  \Pi\setminus\{\gamma\} =\{e^{i\theta_1},\dots, e^{i\theta_{p-1}}\}$
such that
$$0< \theta_1 <\cdots <\theta_m<\pi<\theta_{m+1}<\cdots< \theta_{p-1}<2\pi.$$
As $\Pi\setminus \{\gamma\}$ is $(k-1)$-regular, by Proposition \ref{lower}, $\Pi\setminus \{\gamma\}$ has
$p-1\ge 2(k-1)+1$ elements. Therefore, for every $\alpha\in\Omega$,  the open half plane  $\alpha P$  contains at least $k-1$ elements in $\Pi\setminus\{\gamma\}$ and either $P$ or $- P$ contains at least $k$ elements in $\Pi\setminus\{\gamma\}$. Hence, we have either $m=k-1$ or $k\le m\le p-k$.

Choose $\beta =e^{i\theta}$ where
$$\theta
= \cases{
\max\{\pi+\theta_m,\theta_{p-1}\}/2 & \mbox{ if } $m=k-1$, \cr
\min\{2\pi+\theta_1,\pi+\theta_{m+1}\}/2  & \mbox{ if } $k\le m\le p-k.$
}
$$
Now for every $\alpha \ne \pm 1$, the open half plane $\alpha
P$ contains at least $k-1$ elements of $\Pi\setminus \{\gamma\}$
and either $\gamma$ or $\beta$. Hence, $\alpha P$ contains at least $k$ elements of
$\Pi\cup\{\beta\}$. On the other hand, when $\alpha = \pm 1$, the open half plane
$\alpha P$ contains either $k$ elements of $ \Pi $ or $k-1$  elements of $\Pi $ and $\beta$. Again, $\alpha P$ contains at least $k$ elements of
$\Pi\cup\{\beta\}$. Thus, (a) holds.

\smallskip
Conversely, suppose (a) holds.
If  $-\beta \in \Pi$,
then it is easy to see that the set $\Pi \setminus \{-\beta\}$
is $(k-1)$-regular.
From now, we assume that $-\beta\notin \Pi$.
Without loss of generality, we may assume that $\beta = -1$.
Furthermore, by replacing $\Pi$ with the set $\{\bar \xi: \xi\in \Pi\}$, if necessary, we can assume that
the number of elements in $\Pi \cap P$
is greater or equal to the number of elements in $\Pi \cap (-P)$.
Under this assumption, the upper open half plane must contain at least one
non-antipodal point of $\Pi$.

Let $\gamma$ be the non-antipodal point in $\Pi$ such that
$0<\arg(\gamma) \le \arg(\alpha)$ for all non-antipodal
points $\alpha \in \Pi$. Then $\gamma \in P$.
We show
that $\Pi \setminus \{\gamma\}$ is $(k-1)$-regular.

Take any $\alpha \in \Pi \setminus \{\gamma\}$.
Suppose $\alpha \in \beta P \cup \gamma P$. Then the open half plane $\alpha P$
can contain at most one of points $\beta$ and $\gamma$. As the open half plane
$\alpha P$ contains at least $k$ elements of $\Pi \cup \{\beta\}$,
$\alpha P$ contains at least $k-1$ elements of $\Pi \setminus \{\gamma\}$.
Thus, $\Pi \setminus \{\gamma\}$ is $(k-1)$-regular
if $\Pi\setminus \{\gamma\} \subseteq \beta P \cup \gamma P$.
Now suppose $(\Pi\setminus \{\gamma\}) \setminus (\beta P \cup \gamma P)$
is nonempty and let $\omega_1,\dots,\omega_t$ be the points in this set.
Notice that all of them lie in the upper open half plane $P$.
Therefore, we may assume that $$0 < \arg(\omega_1) < \cdots <
\arg(\omega_t) < \arg(\gamma) < \pi.$$
Also by the choice of $\gamma$,
$\omega_1,\dots,\omega_t$ cannot be non-antipodal points and
hence the points $-\omega_1, \dots,-\omega_t$ are in $\Pi$.
 Clearly, each open half
plane $\omega_j P$ contains at least $k$ elements of $\Pi \cup \{\beta\}$.
Notice that $(w_j P)\setminus P$ contains exactly
$j$ elements of $\Pi \cup \{\beta\}$, namely,
$-\omega_1,\dots,-\omega_{j-1}$ and $\beta$.
Also the set $P\setminus (w_jP)$
contains exactly $j$ elements of $\Pi \cup \{\beta\}$, namely,
$\omega_1,\dots,\omega_j$.
It follows that the half plane $w_jP$ contains the
same number of elements of $P\cup \{\beta\}$
as the upper half plane $P$.
By Proposition \ref{lower},
$\Pi \cup \{\beta\}$ contains at least $2k+2$ elements.
Then by assumption, the upper open half plane $P$ contains at
least $k+1$ elements of $\Pi \cup \{\beta\}$.
Thus, every open half plane $w_j P$ contains at least $k+1$
elements of $\Pi \cup \{\beta\}$ and it
contains at least $k-1$ elements of $\Pi\setminus \{\gamma\}$.
Therefore,
$\Pi\setminus \{\gamma\}$ is a
$(k-1)$-regular set and the assertion follows.
\end{proof}

\smallskip
Applying the above lemma inductively (repeatedly), we have the following.

\medskip
\begin{proposition}\label{prop3}
Let $k > 1$ and $\Pi$
is a $p$ element subset of $\Omega$
containing $s$ pairs of antipodal points,
where $p \ge 3$ and $s \ge 0$.
Suppose $p > 2s$.
For any positive $t \le \min\{k,p-2s,p-1\}$,
the following are equivalent.
\begin{enumerate}
\item[\rm (a)] One can add $t$ points $\beta_1,\dots,\beta_t \notin \Pi$
so that
$\Pi \cup \{\beta_1,\dots,\beta_t\}$ is $k$-regular.
\item[\rm (b)] One can delete $t$ non-antipodal points $\gamma_1,\dots,\gamma_t \in \Pi$ so that
$\Pi \setminus \{\gamma_1,\dots,\gamma_t\}$ is $(k-t)$-regular.
\end{enumerate}
\end{proposition}

\smallskip
\begin{proof}
Clearly, the result holds for $t = 1$ by Proposition \ref{prop2}.
Assume the statement holds for all $\ell < t$.
Suppose $\Pi \cup \{\beta_1,\dots,\beta_t\}$ is $k$-regular.
Let $\Pi_1 = \Pi \cup \{\beta_1\}$.
Then $\Pi_1 \cup \{\beta_2,\dots,\beta_t\}$ is $k$-regular and it follows from the assumption that
one can find $t-1$ non-antipodal points $\gamma_1,\dots,\gamma_{t-1} \in \Pi \cup\{\beta_1\}$
such that $\Pi_1 \setminus \{\gamma_1,\dots,\gamma_{t-1}\}$ is $(k-t+1)$-regular.
If $\beta_1 \notin \{\gamma_1,\dots,\gamma_{t-1}\}$,
by applying Lemma \ref{prop2} to the set $(\Pi \setminus \{\gamma_1,\dots,\gamma_{t-1}\})
\cup \{\beta_1\}$, one can find another non-antipodal point
$\gamma_t \in \Pi_1 \setminus \{\gamma_1,\dots,\gamma_{t-1}\}$
so that $\Pi \setminus \{\gamma_1,\dots,\gamma_{t}\}$ is $(k-t)$-regular.
On the other hand, if $\beta_1$ is one of the $\gamma_j$, say $\beta_1 = \gamma_1$,
then $\Pi \setminus \{\gamma_2,\dots, \gamma_{t-1}\}$ is $(k-t+1)$-regular.
In this case, take an arbitrary element $\gamma_t \in \Pi \setminus \{\gamma_2,\dots, \gamma_{t-1}\}$
and apply Lemma \ref{prop2} to the set
$(\Pi \setminus \{\gamma_2,\dots, \gamma_t\})\cup \{\gamma_t\}$,
one can find another non-antipodal point $\gamma_{t+1}$ so that
the set $\Pi \setminus \{\gamma_2,\dots, \gamma_{t+1}\}$ is $(k-t)$-regular.
Then (b) follows.
The proof of (b) implying (a) can also be done by induction in a similar way.
\end{proof}

\smallskip
Suppose $k < p-s$.
Given a $p$ element subset $\Pi$ of $\Omega$
containing $s$ pairs of antipodal points,
$\beta_1, -\beta_1, \dots, \beta_s,-\beta_s$ with $s > 0$,
which is not $k$-regular,
the set obtained from $\Pi$ by
deleting all $p-2s$ non-antipodal points
is a $(s-1)$-regular set.
On the other hand,
if $\Pi$ does not have any pair of antipodal points, then $k \le p-1$
and one can always delete $k$ elements to form a $0$-regular set.
In both cases, one see that the following
minimum always exist.
\begin{eqnarray*}
q=&& \min \{ t\in \IN:
\Pi \setminus \{ \beta_1, -\beta_1, \dots, \beta_s,-\beta_s \}
\mbox{ has a {\it t}-element }  \cr
&& \hspace{35mm} \mbox{ subset $T$ such that }
\Pi \setminus T \mbox{ is }
(k-t)\mbox{-regular}\, \}.
\end{eqnarray*}
By Proposition \ref{prop3}, one can always add
this minimum number $q$ of points to $\Pi$ to form a $k$-regular set.
Furthermore, this number $q$ is optimal
in the sense that one cannot add fewer than $q$ elements
to do so.

By definition, $q$ is a positive integer
bounded above by $\min\{k,p-2s\}$.
The following proposition gives more information about the minimum value
(\ref{Q2}) in Theorem \ref{main6}.

\medskip\begin{proposition}\label{main6a}
Using the notation in Theorem \ref{main6}.
If $k < p-s$,
then the value $q$ in (\ref{Q2}) exists and satisfies
$$q \le \cases{
k & \hbox{if } (p,s) = (k+1,0) \hbox{ or } $(k+2,1)$, \cr
\min\{k-1,p-2s\} & \hbox{otherwise.}
}$$
Also $q$ is bounded below by $2k+1-p$ or $2k+2-p$,
depending whether $s$ is zero.
Furthermore, $q = 2k+1-p$
if $p \le k+2$ with $s = 0$
and $q = 2k+2-p$ if $p \le k+3$ with $s > 0$.
\end{proposition}

\smallskip
\begin{proof}
The lower bound can be seen easily from Proposition \ref{lower}.
Also the case when $(p,s) = (k+1,0)$ or $(k+2,1)$
has already discussed.
Now we assume that $(p,s) \notin \{(k+1,0), (k+2,1)\}$.
Consider the case when $s \ge 2$.
Take $t = \min\{k-1,p-2s\}$
and delete $t$ non-antipodal elements in $\Pi$.
Then the resulting set is $(s-1)$-regular and hence $(k-t)$-regular
as $k-t = \max\{1,k-p+2s\} \le s-1$.
Thus, $q \le t$.

Next we consider the case when $s = 1$ and $p \ge k+3$.
Let $\{\alpha,\ -\alpha\}$ be the pair of antipodal points in $\Pi$.
Since $\Pi$ is $1$-regular, there are
$\alpha_1\in (\alpha,\ -\alpha)\cap\Pi$ and
$\alpha_2\in (-\alpha,\ \alpha)\cap\Pi$.
Pick another $k-1$ non-antipodal points
$\alpha_3,\dots,\alpha_{k+1}$ in $\Pi$.
The set $\Pi \setminus \{\alpha_3,\dots,\alpha_{k+1}\}$
containing $\{\alpha_1,\alpha_2,\alpha,-\alpha\}$
is $1$-regular. Then $q \le k-1$.

Finally consider the case when $s = 0$ and
$p \ge k+2$. We may assume that $\Pi=\{e^{i\xi_j}:1\le j\le p\}$
with $0=\xi_1<\cdots<\xi_p<2\pi$.
Since $\Pi$ is $1$-regular, we can choose $\ell$ such that
$\xi_\ell=\max\{\xi_j:0<\xi_j<\pi\}$. Then
$S=\{\xi_1,\xi_{\ell},\xi_{\ell+1}\}$ is $1$-regular.
Then any $p-k+1$ subset of $\Pi$ containing $S$ is $1$-regular.
Thus, $q \le k-1$.
\end{proof}

\smallskip
Now we are ready to present the following:

\bf Proof of Theorems \ref{main6} and \ref{main7}. \rm The assertions on
$q$ and $n$ follows by Propositions \ref{ass2} and \ref{prop3}.
For the last assertion in Theorem \ref{main6},
we see that in order to get a $k$-regular set by adding $q$ points to
$\Pi$,  we need to add at least $k-1$ points $e^{i\xi}$, with $0<\xi<\pi$. If $2k+2-p>k-1$,
then $p-3<k$ and we need to add an extra $k-(p-3)$ points $e^{i\xi}$, with $\pi<\xi<2\pi$,
giving a total of $k-1+k-(p-3)=2k+2-p$ points. This proves the equality in (\ref{optimal1}).
The equality in (\ref{optimal2}) now follows readily.
{\vbox{\hrule height0.6pt\hbox{%
   \vrule height1.3ex width0.6pt\hskip0.8ex
   \vrule width0.6pt}\hrule height0.6pt
  }

\smallskip
To close this section, let us illustrate our
results by the following example.

\smallskip
\begin{example}\label{ex3} \rm
Let the polygon $\cP
= \conv\{1,w,w^2,w^3,w^4,w^5,w^6,w^9\}$ with $w = e^{2\pi i/12}$,
see the following.
\begin{center}
\epsfig{file=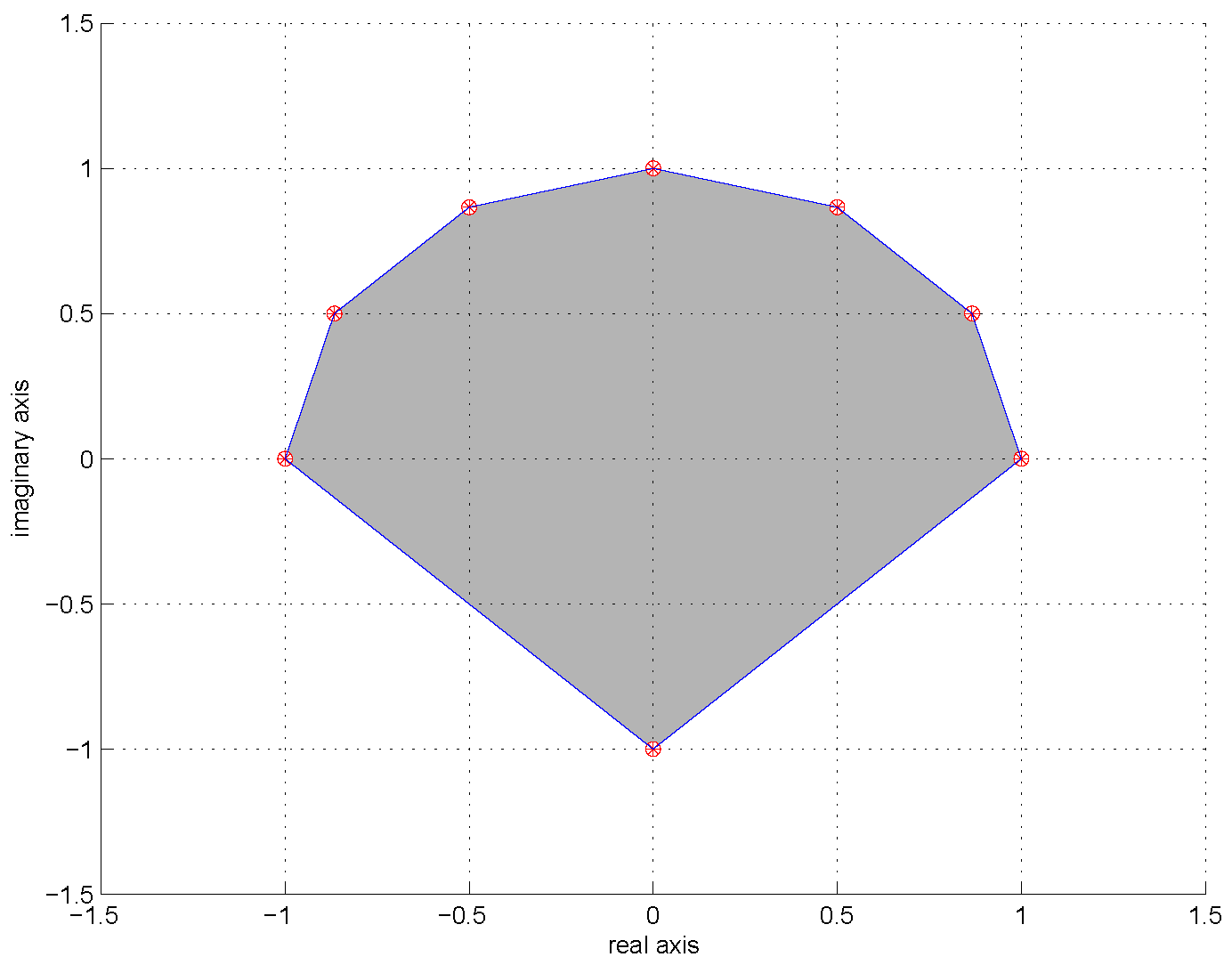, width =1.7in, height=1.7in} \\
The polygon $\cP$
\end{center}
Then $\cP = \bigcap_{j=1}^8 \cH(d_j,\xi_j)$ with
$d_1 = \cdots = d_6 = \cos \frac{\pi}{12}$,
$d_7 = d_8 = \cos \frac{\pi}{4}$, and 
$$\left(\xi_1,\dots,\xi_8\right)
= \left(\frac{\pi}{12}, \frac{3\pi}{12}, \frac{5\pi}{12},
\frac{7\pi}{12}, \frac{9\pi}{12}, \frac{11\pi}{12},
\frac{15\pi}{12}, \frac{21\pi}{12} \right).$$
Thus, 
\vspace{-3mm}
$$\Pi = \{\alpha_1,\dots,\alpha_8\} = \left\{e^\frac{\pi
i}{12}, e^\frac{3\pi i}{12}, e^\frac{5\pi i}{12}, e^\frac{7\pi
i}{12}, e^\frac{9\pi i}{12}, e^\frac{11\pi i}{12}, e^\frac{15\pi
i}{12}, e^\frac{21\pi i}{12} \right\}.$$
In particular,
$\Pi$ has two pairs of antipodal points, namely,
$\left\{e^\frac{3\pi i}{12}, e^\frac{15\pi i}{12} \right\}$ and 
\linebreak
$\left\{e^\frac{9\pi i}{12}, e^\frac{21\pi i}{12} \right\}$,
i.e., $p = 8$ and $s = 2$.
By Theorem \ref{main7} and Proposition \ref{main6a}, for $k \ge 5$,
a $(2k+2)\times (2k+2)$ normal matrix $A$ can be constructed
so that $\Lambda_k(A) = \cP$.

It remains to consider the cases for $k\le 4$.
Clearly, $\Pi$ is $2$-regular.
Thus, a $8\times 8$ normal matrix $A_2$ can be constructed
so that $\Lambda_2(A_2) = \cP$.
However, $\Pi$ is not $k$-regular for $k \ge 3$.

Now we consider the case $k =3$. Clearly, $\Pi\setminus \{ e^\frac{5\pi
i}{12} \}$ is $2$-regular.
Then Theorem \ref{main7} shows that there is a $9\times 9$
normal matrix $A_3$ such that $\Lambda_3(A_3) = \cP$. Indeed, following
the proof of Lemma \ref{prop2}, we see that if $\Pi' =
\{e^\frac{18\pi i}{12} \}$, $\Pi \cup \Pi'$ is
$3$-regular.

Finally, we turn to the case when $k = 4$.
Notice that $\Pi \setminus
\{ e^\frac{5\pi i}{12}, e^\frac{7\pi i}{12} \}$ is $2$-regular. 
Thus, Theorem \ref{main7} shows that there is
a $10\times 10$ normal matrix $A_4$ such that
$\Lambda_4(A_4) = \cP$.

In the following, we display the higher rank numerical ranges of
$A_2,$ $A_3$, and $A_4$. In the figures, the points ``o''
correspond to the vertices of the polygon while the points
``$\ast$'' correspond to the eigenvalues of the normal matrices.

\begin{center}
\epsfig{file=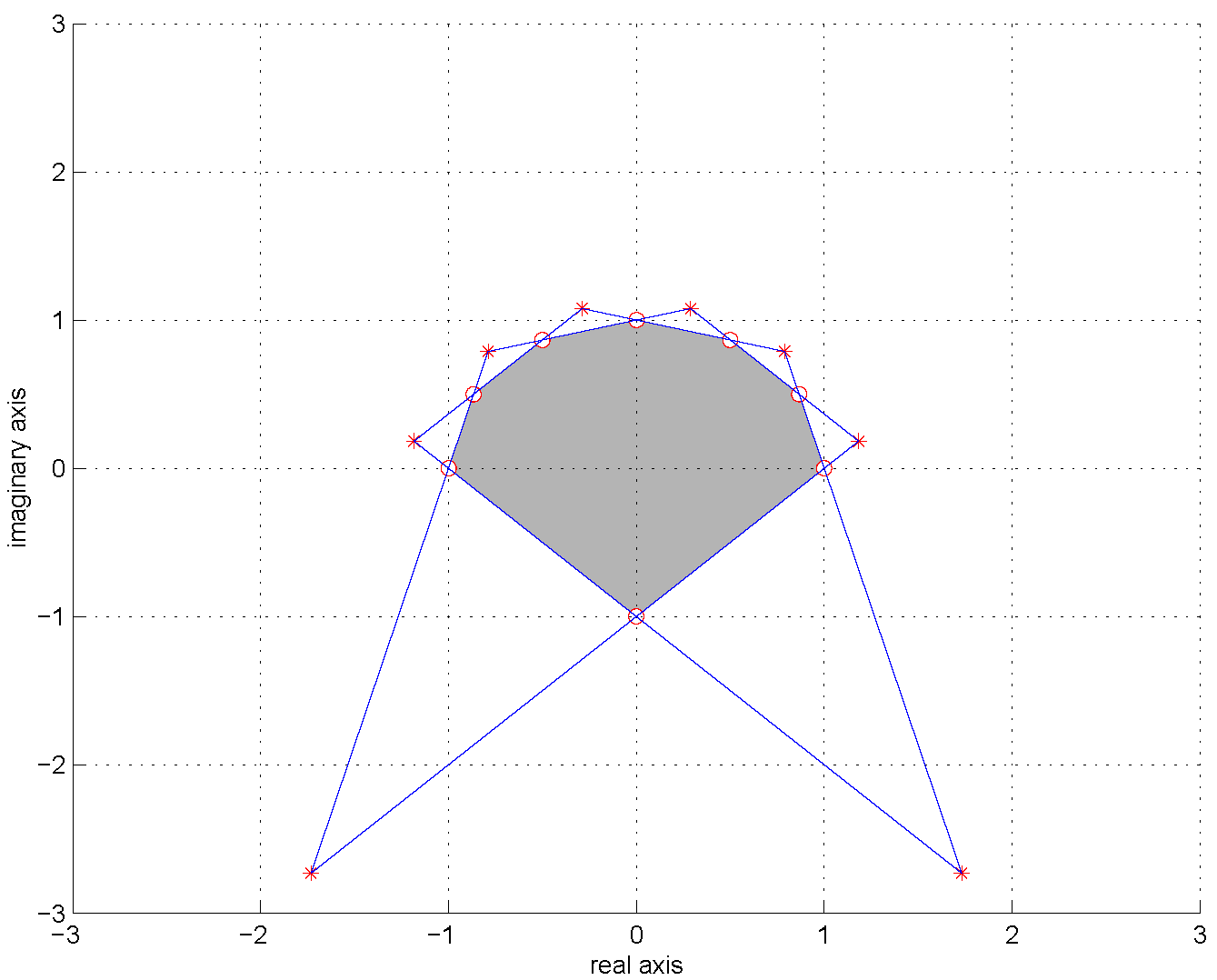, width =2in, height=2in} \qquad
\epsfig{file=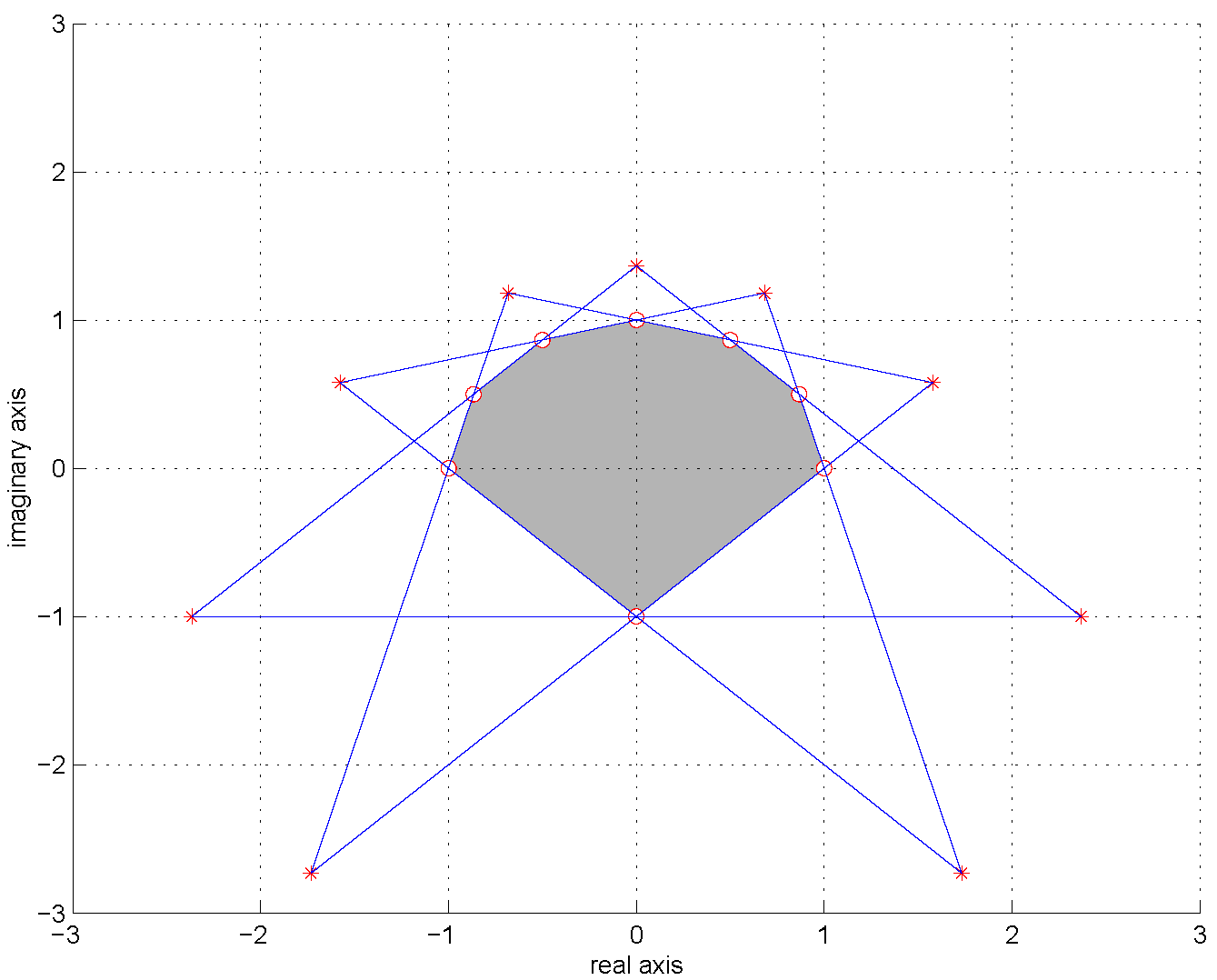, width =2in, height=2in}

$\Lambda_2(A_2) = \cP$ \hspace{3.5cm}
$\Lambda_3(A_3) = \cP$

\epsfig{file=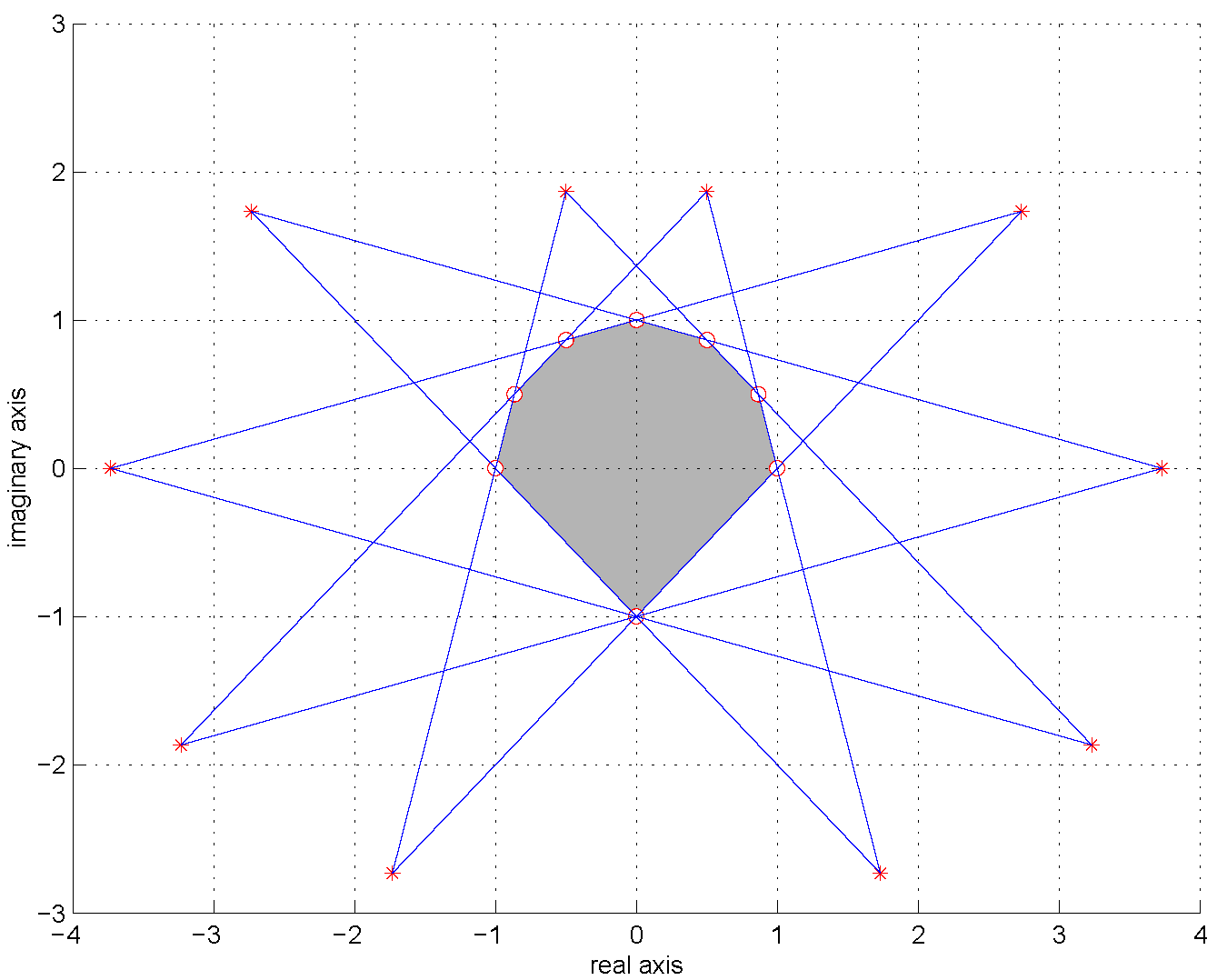, width =2.66in, height=2in}

$\Lambda_4(A_4) = \cP$

\end{center}
\end{example}
\pagebreak

\section{An algorithm}
In this section, we further present a detail procedure 
for constructing the rank-$k$ numerical ranges
of normal matrices based on the discussion in Section 2.

\smallskip
Given a normal matrix $A$ with $m$ distinct eigenvalues
$a_1, \dots, a_m$, one can easily construct $\Lambda_k(A)$ through the
following algorithms.

\medskip
\noindent 
\bf Basic Algorithm \rm 
First construct the set $\cS_0$.
For each ordered pair $(r,s)$ with $r < s$,
count the number of eigenvalues of $A$
(counting multiplicities) in the open planes $H_0(a_r,a_s)$ and $H_0(a_s,a_r)$.
\begin{enumerate}
\item If $H_0(a_r,a_s)$ has at most $n-k-1$ eigenvalues
while $H_0(a_s,a_r)$ has at most $k-1$ eigenvalues,
then collect the index pair $(r,s)$ in $\cS_0$.
\item If $H_0(a_s,a_r)$ has at most $n-k-1$ eigenvalues
while $H_0(a_r,a_s)$ has at most $k-1$ eigenvalues,
then collect the index pair $(s,r)$ in $\cS_0$.
\end{enumerate}

\smallskip
Notice that one can already construct $\Lambda_k(A)$
by determine the intersection of all the half planes
$H(a_r,a_s)$ with $(r,s) \in \cS_0$.
Nevertheless, one can perform the following additional steps 
to simplify the set $\cS_0$ before constructing $\Lambda_k(A)$.

\medskip
\noindent
\bf Modified Algorithm 1 \rm
Suppose in basic algorithm, there is an index pair $(p,q)$ satisfying both (1) and (2),
i.e., both pairs $(p,q)$ and $(q,p)$ are in $\cS_0$.
Then $\Lambda_k(A)$ is a subset of a line segment. 
In this case, $\Lambda_k(A)$ can be constructed as follows.

Set $\hat a_j = (a_j - a_p) / (a_q - a_p)$
and define $\cS_1 = \{(r,s)\in \cS_0: \Im(\hat a_r) \ne \Im(\hat a_s)\}$. If
$\cS_1=\emptyset$, then $\Lambda_k(A)=\emptyset$. Suppose
$\cS_1\ne\emptyset$.
For each $(r,s)\in \cS_1$, compute
$$b_{rs} = \frac{ \Im(\hat a_r)\, \Re(\hat a_s)-\Im(\hat a_s)\, \Re(\hat
a_r) }{\Im(\hat a_r) - \Im(\hat a_s)}.$$
Take
\begin{eqnarray*}
b_1 &=& \max\{b_{rs}: (r,s) \in \cS_1,\ \Im( \hat a_r) \ge 0 \hbox{ and } \Im(\hat a_s) \le 0\}, \cr
b_2 &=& \min\,\{b_{rs}: (r,s) \in \cS_1,\ \Im( \hat a_r) \le 0 \hbox{ and } \Im(\hat a_s) \ge 0\}.
\end{eqnarray*}
Then
$\Lambda_k(A)$ is the line segment in $\IC$ joining the points
$(a_q-a_p) b_1  + a_p$ and $(a_q-a_p) b_2  + a_p$ 
if $b_1 \le b_2$; otherwise, $\Lambda_k(A) = \emptyset$.

\medskip
\noindent
\bf Modified Algorithm 2 \rm
Assume the situation mentioned in modified algorithm 1 does not hold.
Check if the set $\cS_0$ satisfy the following.
\begin{eqnarray}\label{cond3}
&&\hbox{There are } (r_1,s_1), \dots, (r_\ell,s_\ell) \in \cS_0 \hbox{ with  $\ell \ge 3$  such that} \cr
&&\hspace{3cm}
\{r_1,s_1\} \cap \{r_2,s_2\}\cap \cdots \cap \{r_\ell, s_\ell\} = \{t\}
\quad \hbox{for some $1\le t\le m$.}
\end{eqnarray}
If yes, define
$$\theta_j  = \cases{
\arg(a_{s_j} - t) & \hbox{if } $r_j = t$,  \cr
\arg(t - a_{r_j}) & \hbox{if } $s_j = t$. 
}$$
Relabel the indices so that $0 \le \theta_1 \le \cdots \le\theta_\ell < 2\pi$.
Consider the following three cases.
\begin{enumerate}
\item If $\theta_\ell - \theta_1 < \pi$, 
remove the all pairs $(r_j,s_j)$ in $\cS_0$ for $j \ne 1,\ell$.
Then check again whether the modified set still satisfies (\ref{cond3}).
\item If $\theta_{k+1} - \theta_k > \pi$ for some $k$,
remove the all pairs $(r_j,s_j)$ in $\cS_0$ for $j \ne k,k+1$.
Then check again whether the modified set still satisfies (\ref{cond3})
\item If the above two items are not satisfied, 
then $\Lambda_k(A)$ is either the empty set or the singleton set $\{a_t\}$.
In this case, check whether $a_t$ lies in $H(a_r,a_s)$ for all $(r,s) \in \cS_0$.
If yes, $\Lambda_k(A)$ is the singleton set; otherwise it is the empty set.
\end{enumerate}
Finally, if the modified set $\cS_0$ does not satisfy (\ref{cond3}), 
then one can construct $\Lambda_k(A)$
by determine the intersection of all the half planes
$H(a_r,a_s)$ with $(r,s)$ in the modified set $\cS_0$.

\medskip
\noindent
{\bf Acknowledgment}

This research began at the 2008 IMA PI Summer Program for Graduate
Students, where the second author is a lecturer and the third
author is a co-organizer.  The support of IMA and NSF for the
program is graciously acknowledged. The hospitality of the
colleagues at Iowa State University is deeply appreciated.
The authors would also like to thank the referees for some
helpful comments.

\end{document}